\documentclass[hidelinks,onefignum,onetabnum]{siamart251216}
\usepackage{graphicx} %
\usepackage[]{mdframed}
\usepackage{amsfonts}
\usepackage{amsmath,amstext,amsbsy,amssymb,mathdots}
\usepackage{mathrsfs}
\usepackage{appendix}
\usepackage{longtable}
\usepackage{array}
\usepackage{booktabs}
\usepackage{ragged2e}
\usepackage{tabularx}
\usepackage{siunitx}
\usepackage{algorithm}
\usepackage{algpseudocode}
\makeatletter
\let\algx@For\For       \let\algx@EndFor\EndFor   \let\algx@If\If
\let\algx@Else\Else     \let\algx@While\While     \let\algx@EndWhile\EndWhile
\let\algx@Repeat\Repeat \let\algx@Until\Until     \let\algx@Return\Return
\makeatother
\usepackage[algo2e,boxed,linesnumbered]{algorithm2e}
\usepackage{etoolbox}
\makeatletter
\AtBeginEnvironment{algorithmic}{%
  \let\For\algx@For       \let\EndFor\algx@EndFor \let\If\algx@If
  \let\Else\algx@Else     \let\While\algx@While   \let\EndWhile\algx@EndWhile
  \let\Repeat\algx@Repeat \let\Until\algx@Until   \let\Return\algx@Return}
\makeatother

\SetKwFunction{LeftStep}{LeftStep}
\SetKwFunction{RightStep}{RightStep}
\SetKwFunction{AlternatingStep}{AlternatingStep}
\SetKwProg{Fn}{Function}{:}{end}
\SetKwProg{Proc}{Procedure}{:}{end}
\renewcommand{\arraystretch}{1.3}
\theoremstyle{nonumberplain}
\theoremheaderfont{\normalfont\itshape}
\theorembodyfont{\normalfont}
\theoremseparator{.}
\theoremsymbol{}
\newtheorem{remark}{Remark}
\newsiamremark{example}{Example}

\crefname{subsubsection}{subsection}{subsections}
\Crefname{subsubsection}{Subsection}{Subsections}

\newcommand{\sgn}{{\rm sgn}}

\newcommand{ \rank}{{\rm rank}}

\newcommand{\C}{\mathbb{C}}
\newcommand{\RR}{\mathcal{R}}
\newcommand{\G}{\mathcal{G}}
\newcommand{\Gt}{\tilde{\mathcal{G}}}

\headers{A fast solver for lopsided Sylvester equations}{L. Batakci and H. Wilber}

\title{A fast iterative solver for lopsided Sylvester matrix equations\thanks{Submitted to the editors \today.
\funding{This work was supported by NSF grant DMS-2410045.}}}
\author{Levent Batakci\thanks{Department of Applied Mathematics, University of Washington, Seattle, WA 98195 (\email{levent@uw.edu}, \email{hdw27@uw.edu}).}
\and Heather Wilber\footnotemark[2]}

\begin{document}
\maketitle

\begin{abstract}
A generalization of the factored alternating direction implicit algorithm is introduced that makes use of type $(p, q)$ Zolotarev rational functions to apply unequal numbers of left and right shifts. Convergence properties of the algorithm are discussed, as well as a method for computing shift parameters. Several applications are highlighted where the generalized method is useful, including in the solving of a collection of shifted linear systems, for discretizations that involve left and right operators that are not equally amenable to inversion, and for Sylvester matrix equations with differently-sized, differently-structured, or differently-conditioned coefficient matrices. 
\end{abstract}


\section{Introduction}

Fast and memory-efficient methods for solving Sylvester matrix equations of the form
\begin{equation} \label{eq:sylv} AX - XB = C,
\end{equation}
where $X, C \in \mathbb{C}^{m \times n}$, $A \in \mathbb{C}^{m \times m}$ and
$B \in \mathbb{C}^{n \times n}$, are centrally important in many applications in
scientific computing. This includes dynamical systems and control~\cite{antoulas2005approximation,simoncini2016computational}, signal and image processing~\cite{calvetti1996application}, computational methods for solving partial differential equations (PDEs)~\cite{FortunatoTownsend2020,PalittaSimoncini2016, PowellSilvesterSimoncini2017}, eigenvalue assignment problems~\cite{BhattacharyyaDeSouza1982,CavinBhattacharyya1983}, and the solution of large-scale algebraic Riccati equations by Newton--Kleinman iteration~\cite{BennerSaak2013}. 
In these settings, $m$ and $n$ can be extremely large, and the goal is to efficiently compute an approximation to the unknown matrix $X$. 
In many of these applications, the \textit{coefficient matrices} $A$ and $B$ are data-sparse, and $C$ is low rank~\cite{Baur2008,Grasedyck2004}. 
The unknown matrix $X$ is typically dense, but often is well approximated by low-rank matrices~\cite{BeckermannTownsend2017, Grasedyck2004,Penzl2000}. 
Direct solvers for~\eqref{eq:sylv} require $\mathcal{O}(mn)$ storage complexity and $\mathcal{O}(m^3 + n^3)$ computational complexity~\cite{BartelsStewart1972}. 
In the large-scale setting, these problems are typically approached via iterative methods that construct approximations to $X$ in low-rank form.

In this paper, we consider settings where the computational cost $T_A$ of applying a shifted invert with $A$ is substantially greater than $T_B$, or vice versa. Analogously, there are settings where memory, access, or conditioning issues make the inversion of $B$ more amenable than that of $A$.  We call such equations \textit{lopsided} Sylvester matrix equations. Examples include collections of shifted linear systems, which also arise when evaluating functions of matrices such as $A^{-1/2}b$ (\cref{sec:ex:shifted}), as well as various PDE discretizations that result in coefficient matrices with drastically different sizes, structures, and properties.  %

A natural idea is to design iterative methods that limit or avoid inverts involving $A$ while still benefiting from inverts involving $B$. We refer to such methods as \textit{lopsided iterative methods}. A memory-dense lopsided version of the alternating direction implicit (ADI) method was introduced by Levin and Saff in~\cite{levin1994optimal} (the g-ADI method). We expand on this idea and develop a memory-efficient lopsided iterative solver called g-fADI that constructs approximate solutions to~\eqref{eq:sylv} in low-rank form.
In the nonlopsided setting, g-fADI is equivalent to the factored ADI method~\cite{benner2009adi, li2002low}. The generalization turns out to be nontrivial; an efficient implementation benefits from the splitting-based framework introduced in~\cite{schulze2024unifying}.
Various settings for which g-fADI is a good fit are shown in \cref{tab:usecases}.

To analyze the convergence properties of g-fADI, we expand on Levin and Saff's analysis of the so-called \textit{generalized Zolotarev rational approximation problem} that describes its error behavior~\cite{levin1994optimal}, and we introduce a numerical method for computing generalized Zolotarev rationals based on the approach in~\cite{trefethen2024computation}. We also describe a collection of applications where the lopsided algorithms are useful, including settings where $A$ is large and $B$ is small, as well as settings where $A$ is structured and $B$ is not. In particular, we show that systems of shifted linear equations of the form $(A+s_jI)x_j = b$ can be recast as lopsided Sylvester matrix equations that are often efficiently solved with g-fADI using only products with $A$. %

\begin{table}[t]
\caption{Settings where lopsided methods, such as g-fADI, can be useful.}
\label{tab:usecases}
\centering\small
\renewcommand{\arraystretch}{1.15}
\begin{tabularx}{\textwidth}{@{}>{\raggedright\arraybackslash}Xl@{}}
\toprule
Use case & See \\
\midrule
Many shifted linear systems $(A+s_jI)x_j=b$, e.g., frequency sweeps & \cref{sec:ex:shifted} \\
Evaluation of $f(A)b$ via rational approximation & \cref{sec:ex:shifted} \\
Large spatial operator coupled to a small one in time or a parameter & \cref{sec:ex:pde} \\
Diffusions that switch among the states of a Markov chain & \cref{sec:ex:switch} \\
Nonsymmetric operators, where Krylov bases are costly to store & \cref{sec:ex:switch,sec:ex:sweep} \\
Small coefficient matrix that is  far from normal & \cref{sec:ex:pde} \\
Only a few outputs $\Gamma X$ are needed (e.g., sensor values) & \cref{sec:ex:switch,sec:ex:sweep} \\
Storing the iterates with other methods exceeds available memory & \cref{sec:ex:pde} \\
\bottomrule
\end{tabularx}
\end{table}

\subsection{A model problem}\label{sec:intro:teaser}
To illustrate the kinds of applications where lopsided Sylvester matrix equations arise, consider  $AX-XB=UV^*$ with $B=-L^T$, where $A$ is the second-order centered 7-point discretization of $-\Delta+w\cdot\nabla$ on $N^3$ interior points of the unit cube, with a constant drift $w$ of length 5 in the direction $(1,0.6,0.3)$~\eqref{eq:tbK}, and $L\in\mathbb{R}^{200\times200}$ is the shifted and scaled second-order finite difference discretization~\eqref{eq:tbL} of a one-dimensional Laplace operator. The right-hand side is random and of rank two, and the target accuracy is $10^{-8}$. Thus $A$ is large, sparse, and nonsymmetric, while $B$ is small, symmetric, and tridiagonal.

This equation describes the steady state of a population that diffuses through a homogeneous three-dimensional medium, drifts with the constant velocity $w$, and is structured by a phenotypic trait $t\in(0,1)$, in which mutations act as diffusion in $t$. Its density satisfies
\[
  -\Delta u+w\cdot\nabla u-\epsilon\,\partial_t^2u+u=f,\qquad (x,t)\in(0,1)^3\times(0,1),
\]
with homogeneous Dirichlet conditions and $\epsilon\approx0.062$. This is a linear, constant-coefficient version of models of populations structured by space and a phenotypic trait~\cite{AlfaroBerestyckiRaoul2017, AlfaroCovilleRaoul2013, Perthame2007}, with a constant death rate in place of growth and competition.  More generally, equations of this form arise from tensor-product discretizations of problems on a product domain $\Omega\times(0,Y)$, where the coordinate over $(0, Y)$ needs far fewer unknowns than the three spatial ones. The extension approach to fractional diffusion leads to the same structure~\cite{Banjai2019tensor, CaffarelliSilvestre2007, NochettoOtarolaSalgado2015}, as do additional examples we describe in \cref{sec:examples}. 

\begin{figure}[t!]
\centering
\includegraphics{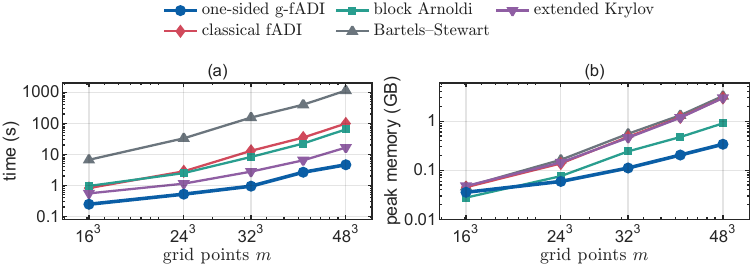}
\caption{(a) Run time and (b) peak memory for $AX+XL^T=UV^T$, with $A$ the discretized $-\Delta+w\cdot\nabla$ on $m=N^3$ grid points and tridiagonal $L$ of size 200, to a relative residual of $10^{-8}$, against $m$. One-sided g-fADI uses only products with $A$. All setup costs are included in the timings.}
\label{fig:teaser}
\end{figure}

The g-fADI method we introduce here supplies a low-rank approximation to the solution and never applies solves involving $A$ to do so. Instead, each step uses one product with $A$ and one tridiagonal solve with $B$. We also benefit from a specialized compression method we introduce called ``small side compression'' (see \cref{alg:onesided}). \Cref{fig:teaser} shows the results against a collection of other methods and illustrates the potential advantages in both speed and memory that can be gained with our lopsided method. While g-fADI takes many more iterations than the other methods,  each step is cheap, and the basis for the smaller low-rank factor in the solution never has more than 59 of its 200 possible columns. Moreover, the a priori bound we use to determine the number of g-fADI iterations is conservative. It results in residuals on the order of  $10^{-12}$, well below the target accuracy of $10^{-8}$. 

\subsection{Generalized Zolotarev rational functions} 
\begin{figure}[t!]
    \centering
    \includegraphics{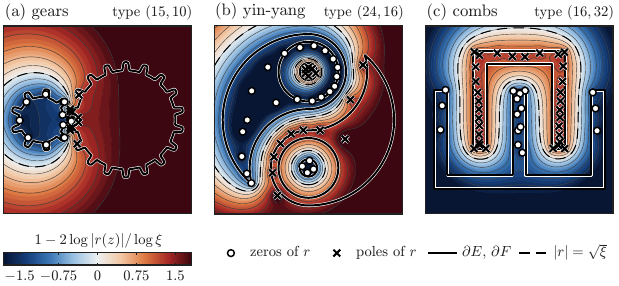}
    \caption{Plots of $1-2\log|r(z)|/\log \xi_{p,q}(E,F)$, where $r$ is a type $(p,q)$ generalized Zolotarev rational, are shown on a log-scale in the complex plane for different choices of $E$ and $F$.  When $|r(z)|$ is small, the plotted function is $<0$ (blue), and when it is large, $>0$ (red).  The poles (x marks) and zeros (white circles) of the rational are parameters required to run the g-fADI algorithm. }
    \label{fig:zolo_gallery}
\end{figure}

We call a rational function in irreducible form with $p$ zeros and $q$ poles (counting multiplicities) a type-$(p,q)$ rational function. Critical in our analysis are type-$(p,q)$ Zolotarev rational functions, which are sometimes referred to as \textit{generalized} Zolotarev rationals~\cite{LeBaillyThiran2000}. 
These functions arise in the error equation for g-fADI, and they also appear in a number of other applications in computational math~\cite{ BeckermannTownsend2019,GuttelPolizzi2015,KressnerSusnjara2017,NakatsukasaFreund2016}. 
Much work has focused on theoretical properties of the Zolotarev rationals~\cite{Gonchar1969, istace1995third, LeBaillyThiran2000, Lebedev1977, Starke1992}, usually in the type-$(p,p)$ or type-$(p\!-\!1,p)$ case. 
For a limited collection of domains, analytical arguments or heuristics can be used to computationally construct the rationals~\cite{DruskinKnizhnermanZaslavsky2009,Wachspress2013}. 
Recently, a numerical method for constructing type-$(p,p)$ Zolotarev rationals via the adaptive Antoulas--Anderson (AAA) algorithm was introduced~\cite{trefethen2024computation}, along with a method based on the Loewner framework~\cite{poussotvassal2025zolotarev}.
We expand on the AAA-based method to compute type-$(p,q)$ Zolotarev rationals (see \cref{fig:zolo_gallery}). 
The poles and zeros of these rational functions supply critical parameters for running the g-fADI and g-ADI algorithms.

We state Zolotarev's central rational approximation problem (typically referred to as \textit{Zolotarev's third problem} below). Except in a handful of important cases~\cite{Akhiezer1990,Gonchar1969,Starke1992,Zolotarev1877}, explicit solutions to these problems are unknown, even in the type-$(p,p)$ case. 
Computational solutions are therefore essential.

\vspace{.5cm}

\begin{mdframed}
\noindent \textbf{(Z3) Zolotarev's third problem.} Given disjoint sets $E$ and $F$ in the complex plane and fixed nonnegative integers $p,q$, find a rational function $r$ of type $(p,q)$ such that 
\begin{equation}\label{eq:Zolotarev_ratio}
    \xi(r) = \frac{\sup_{z \in E} |r(z)|}{\inf_{z\in F}|r(z)|}
\end{equation}
is minimized.
The minimum 
\begin{equation}\label{eq:Z3_number}
    \xi_{p,q}(E,F):= \min_{r \in \mathcal{R}_{p,q}}\xi(r)
\end{equation}
is called the $(p,q)$th Zolotarev number associated with $(E,F)$.
\end{mdframed}
\vspace{.5cm}
The essential property of a Zolotarev rational function, illustrated in \cref{fig:zolo_gallery}, is that it becomes as small as possible on one set ($E$), while remaining relatively large (e.g., of magnitude $> 1$) on the other set ($F$).  In this way,  Zolotarev's third problem is the rational analogue of the Chebyshev minimax polynomial approximation problem underlying GMRES and CG~\cite{driscoll1998}.  Our computational method handles challenging domains (e.g., multiply-connected regions, interlacing narrow inlets,  corners and cusps) and only needs boundary information for the sets $E$ and $F$ in order to be applied. 

The properties of $\xi_{p,q}$ are essential to the analysis of the g-fADI and g-ADI methods.
If $X_{p+q}$ is the approximation produced by g-fADI after $q$ and $p$ shifted inverts of $A$ and $B$, respectively, with shifts chosen as the zeros and poles of the type-$(p,q)$ Zolotarev rational, then we have the error bound
\[
    \|X_{p+q}-X\|_2 \le 4\xi_{p,q}\|X\|_2.
\]
The factor 4 can be reduced to 2 or 1 if one or both of $A$ and $B$ is normal, respectively.
Avoiding expensive steps is the main goal of lopsided solvers, and the behavior of $\xi_{p,q}$ supplies a method for analyzing exactly how many expensive inversion steps we can avoid without suffering from problematic convergence issues.
In many cases, we find that expensive inversion steps can be avoided altogether.

\vspace{.5cm}

\subsection{Organization}
The remainder of this paper is organized as follows. In \cref{sec:adi}, we derive the g-fADI algorithm. In \cref{sec:theory}, we describe theoretical results for generalized Zolotarev rationals and a method for computing them numerically. \Cref{sec:examples} presents numerical examples and applications. MATLAB code for g-fADI and for all of our experiments is available on GitHub~\cite{gfadiCode}.

%%%%%%%%%%%%%%%%%%%%%%%%%%%%%%%%%%%%%%%%%%%%%%%%%%%%%%%%%%%%%%%%%%%%%%%%%%%%%%%%%%%%%%%%%%%%%%%%%%%%%%%%%%%%%
%------------------------ADI AND ZOLOTAREV------------------------------------------------------------------%
%%%%%%%%%%%%%%%%%%%%%%%%%%%%%%%%%%%%%%%%%%%%%%%%%%%%%%%%%%%%%%%%%%%%%%%%%%%%%%%%%%%%%%%%%%%%%%%%%%%%%%%%%%%%%
\section{A low-rank generalized ADI method}\label{sec:adi}

Many numerical methods have been developed for solving \eqref{eq:sylv}, and a detailed review can be found in the survey~\cite{simoncini2016computational}. Before discussing the g-ADI method and expanding it into a low-rank version, we supply a brief overview of existing methods that helps identify where ADI-based methods are useful. 

\subsection{Direct vs. splitting-based vs. projection-based methods} The primary direct method for solving~\eqref{eq:sylv} is the Bartels--Stewart algorithm, having $O(mn)$ storage complexity and $O(m^3 + n^3)$ computational complexity~\cite{BartelsStewart1972}. 
An improvement offering a constant-factor speedup was developed in~\cite{golub1979hessenberg}, but both storage and arithmetic costs remain prohibitive for large problems.

For large problems and/or settings where $A$ and $B$ are banded, sparse, or otherwise structured, iterative methods that exploit often-present  numerical low-rank structure in the solution can be used instead.
The two main types of iterative methods are \textit{projection-based} and \textit{splitting-based}.
Projection-based methods iteratively search for $X$ in a low-dimensional subspace, typically a rational Krylov subspace, computing it as the solution to a smaller, projected Sylvester matrix equation~\cite{CasulliRobol2024, JaimoukhaKasenally1994, KressnerLundMasseiPalitta2021, Saad1990, simoncini2007}.
Their efficiency depends on the choice of subspace, which can be updated adaptively at each step using Ritz values~\cite{druskin2011adaptive} or tangential interpolation~\cite{DruskinSimonciniZaslavsky2014}.

In contrast, splitting methods produce approximate solutions through a sequence of shifted solves and shifted products.
The alternating-direction-implicit (ADI) iteration is the prototype and was introduced in the context of numerically solving the 2D Poisson equation~\cite{peaceman1955numerical,Wachspress2013}. 
When the right-hand side $C$ is low rank, there is a factored version (fADI) which builds $X$ directly in low-rank form~\cite{li2002low, Penzl2000cyclic}. 
Convergence depends critically on the quality of the shifts~\cite{BennerKurschnerSaak2013,Wachspress2013}, and there is a wealth of heuristics available for choosing them~\cite{Penzl2000cyclic, sabino2007solution, Starke1992}. We discuss shift selection for our method in \cref{sec:computing shifts}.

These two classes of iterative methods are related.
In both, an error bound can be obtained from Zolotarev's third problem on the spectra or fields of values of $A$ and $B$~\cite{beckermann2011error,benner2009adi,druskin2011analysis},
with the methods coinciding for particular choices of parameters~\cite{druskin2011analysis, flagg2013adi}.
Projection-based methods have higher computational complexity, owing to the cost of subspace orthogonalization (which grows at every step), but can  achieve much better convergence rates.
ADI is thus attractive when good shifts, which guarantee fast convergence, can be computed, which in turn requires some knowledge about the fields-of-values or spectral sets of the coefficient matrices $A$ and $B$. This is often reasonable in settings where $A$ and $B$ come from the discretization of operators~\cite{FortunatoTownsend2020, shi2021compressibility, wilber2017computing}, or in settings where  precomputation devoted to probing and/or estimating the fields-of-values is reasonable. We supply several examples in \cref{sec:examples}.  By computing or estimating $\xi_{p,q}(E, F)$, where $E$ and $F$ enclose the fields of values of $A$ and $B$, the convergence behavior of ADI-based methods is guaranteed a priori and the number of steps can be assigned without monitoring a residual. Moreover, every step has the same cost, and expensive orthogonalization and projection routines are not needed. Even when a projection-based method is preferred (e.g., because the spectral content of $A$ and $B$ can only be poorly estimated), the convergence theory of ADI remains useful, since the fastest converging projection methods, based on rational Krylov subspace iterations, have errors controlled by the ADI error~\cite{beckermann2011error, druskin2011analysis, druskin2011adaptive}.  Hybrid methods that combine ADI-like and projection-like steps are also available~\cite{benner2009adi}. 

\subsection{Generalized ADI and error bounds}
The standard ADI iteration proceeds via the alternating updates,
    \begin{subequations}\label{eq:ADI}
        \begin{align}
        (A-\beta_{i+1}I)X_{i+1/2} &= X_i(B-\beta_{i+1}I) + C, \label{eq:ADI-left} \\
        X_{i+1}(B-\alpha_{i+1}I) &= (A-\alpha_{i+1}I)X_{i+1/2} -C, \label{eq:ADI-right}
        \end{align}
    \end{subequations}
where the $\alpha_i \in \C$ and $\beta_i \in \C$ are the so-called \textit{shift parameters}.
We refer to \eqref{eq:ADI-left} and \eqref{eq:ADI-right} as left steps and right steps, respectively, as \eqref{eq:ADI-left} involves a solve on the left side of the iterate, and \eqref{eq:ADI-right} involves a solve from the right side.
Theorem 4.1 of~\cite{li2002low} shows that the final approximation is invariant to permutation of the shifts.
This property was recently generalized to fully-commuting splitting schemes in Proposition 2.4 of~\cite{schulze2024unifying}.
An immediate consequence is that in infinite precision, the order of the steps does not matter: the final approximation will be the same if the iterations are completed in alternating order as shown above, or if they are completed with all the left steps taken first, followed by all the right steps.
This is not the entire story. \Cref{fig:shift_ordering} shows that the computation can be numerically sensitive to ordering of the shifts; this phenomenon is also discussed in~\cite{LebedevFinogenov1971, sabino2007solution}. 

Levin and Saff proposed a generalization of ADI (g-ADI) in~\cite{levin1994optimal}, where the numbers of left and right steps can differ.
Let the sets of parameters be denoted as $\{\alpha_i\}_1^{p}$ and $\{\beta_i\}_1^{q}$.
After a total of $p$  right steps and $q$  left steps, the g-ADI error is given by
    \begin{equation}\label{eq:generalized-ADI-error}
        X_{p+q}-X = r_{p,q}(A)(X_0-X)r_{p,q}^{-1}(B), 
        ~~\text{ where } 
        r_{p,q}(z) = \prod_{i=1}^p (z-\alpha_i) \Big/ \prod_{i=1}^q(z-\beta_i),
    \end{equation}
from which we obtain the error bound
\[
    \frac{\|X_{p+q}-X\|}{\|X_0-X\|} \le \|r_{p,q}(A)\|~\|r_{p,q}^{-1}(B)\|.
\]
This bound can be recast in terms of $k$-spectral sets\footnote{A set $\Omega \subseteq \C$ is said to be a $k$-spectral set of a matrix $A$ if (i) $\Lambda(A) \subseteq \Omega$, where $\Lambda(A)$ is the spectrum of $A$, and (ii)~$\|r(A)\|~\le~k\max_{z \in \Omega}|r(z)|$ for all rational functions $r(z)$ analytic on $\Omega$.} of $A$ and $B$.
With the recent resolution of Crouzeix's conjecture~\cite{jin2026numerical,lorist2026}, we have
\begin{equation}\label{eq:generalized-ADI-error-bound-nonnorm}
        \frac{\|X_{p+q}-X\|}{\|X_0-X\|}
        \le 4
        \frac{
            \max_{z \in W(A)}|r_{p,q}(z)|}
            {\min_{z \in W(B)}|r_{p,q}(z)|},
    \end{equation}
where $W(A)$ and $W(B)$ are the fields of values of $A$ and $B$, respectively.
If $A$ and $B$ are normal matrices, the bound becomes
    \begin{equation}\label{eq:generalized-ADI-error-bound}
        \frac{\|X_{p+q}-X\|}{\|X_0-X\|}
        \le  
        \frac{
            \max_{\lambda \in \Lambda(A)}|r_{p,q}(\lambda)|}
            {\min_{\lambda \in \Lambda(B)}|r_{p,q}(\lambda)|},
    \end{equation}
where $\Lambda(A)$ and $\Lambda(B)$ are the spectra of $A$ and $B$, respectively.
The problem of choosing good shift parameters is equivalent to choosing a rational function that minimizes the bound in~\eqref{eq:generalized-ADI-error-bound-nonnorm}, i.e., solves the type-$(p,q)$ Zolotarev ratio problem on a $k$-spectral set of $A$ and $B$.
Particularly, the poles and zeros of an infimizing rational supply  the $\beta_j$'s and $\alpha_j$'s, respectively.
We discuss computational methods for solving this problem and thus retrieving shift parameters  in \cref{sec:computing shifts}. 
For the remainder of this section, we assume shift parameters are given and focus on how to develop a low-rank version of g-ADI.
\subsection{Standard factored ADI}
We begin with a review of factored alternating-direction-implicit (fADI) in the setting where $p=q$.
Our method generalizes on this construction for $p \neq q$.
A useful property of ADI is that the rank of the iterates grows slowly if $C$ is low rank. 
Indeed, collapsing \eqref{eq:ADI-left} and \eqref{eq:ADI-right} shows that $\rank(X_{j+1}) \le \rank(X_j) + \rank(C)$, where $\rho=\rank(C)$.
A more careful analysis establishes the stronger result,
\begin{equation}\label{eq:rank-update}
    X_{j+1}=X_{j} + \Delta X_{j+1}, ~~~\text{where }~ \rank(\Delta X_j) = \rho,
\end{equation}
showing that each iteration of ADI is a rank-$\rho$ update~\cite{li2002low, schulze2024unifying}.
Starting with a low-rank initial guess, commonly $X_0=0$, the ADI algorithm builds low-rank approximations to $X$.
This motivates constructing  $X_k$ directly via its low-rank factorizations. In the case where $p=q$,  this algorithm is known as factored ADI (fADI).
It was first developed for Lyapunov matrix equations with cyclic shifts in~\cite{Penzl2000cyclic}, then generalized to noncyclic shifts and optimized in~\cite{li2002low}, and finally generalized to Sylvester matrix equations in~\cite{benner2009adi}.
Schulze and Saak extended the fADI algorithm to allow for $X_0\ne 0$ and demonstrated that several applications benefit greatly from a low-rank, nonzero initialization~\cite{schulze2024unifying}.
We assume $X_0=0$ for the rest of the paper, but our method can readily be extended to allow for nonzero initializations.
 
Let $C=UV^*$ with $U \in \C^{m \times \rho}$ and $V \in \C^{n \times \rho}$ so that $\rank(C)=\rho$. 
The $k$th iterate is given by $X_k = \sum_{1 \le j \le k} \Delta X_j$, where the updates admit rank-$\rho$ factorizations $\Delta X_j = \gamma_jZ_jY_j^*$.
The fADI algorithm~\cite{BennerKurschner2014} computes these factorizations as 
\begin{equation}\label{eq:standard-fADI}
\gamma_j = (\beta_j-\alpha_j)
~\text{ and }~
    \begin{cases}
        Z_1=(A-\beta_{1}I)^{-1}U,\\
        Z_{j+1} = Z_j+(\beta_{j+1}-\alpha_j)(A-\beta_{j+1}I)^{-1}Z_j,\\
        Y_1=(B^*-\bar\alpha_{1}I)^{-1}V,\\
        Y_{j+1} = Y_j+(\bar\alpha_{j+1}-\bar\beta_j)(B^*-\bar\alpha_{j+1}I)^{-1}Y_j.\\
    \end{cases}
\end{equation}
The fADI updates have computational complexity $O(\rho(T_A+T_B))$ and storage complexity $O(\rho(m+n))$.
In contrast, iterations of standard ADI~\eqref{eq:ADI} have computational complexity $O(nT_A + mT_B)$ and storage complexity $O(mn)$.

\subsection{Generalized factored ADI}\label{subsec:gfadi}
The ideas behind fADI can be applied to the generalized ADI 
    framework and offer similar benefits.
In this section, we derive update rules that enable the
    construction of the iterates $X_j$ in low-rank form, for
    arbitrary ordering of left and right steps.
Our derivation relies on the framework for ADI-like splitting methods proposed in~\cite{schulze2024unifying}.
Since we will be using an unequal number of left and right steps, we drop the half-index notation present in standard ADI.
Each step (left or right) increases the index of the iterate by one.

\begin{figure}[t!]
    \centering
    \includegraphics{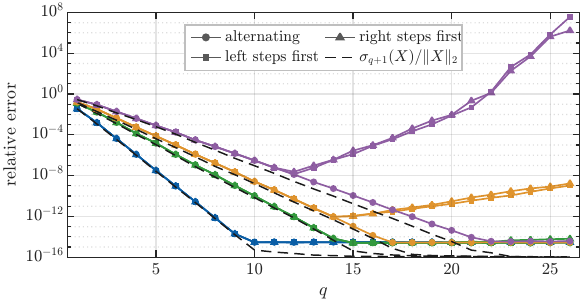}
    \caption{Four g-fADI experiments with $A,B,C \in \mathbb{R}^{300\times300}$ and $\rank(C)=1$.
    The dashed line shows the singular values of $X$; solid lines show the error in the
    rank-$q$ approximants produced by g-fADI with the type-$(q,q)$ Zolotarev shifts for $[-1,-d]$ and $[d,1]$ under three step orderings: alternating, all left steps first, and all right steps first. In each experiment, $\Lambda(B),-\Lambda(A) \subseteq[d,1]$ with $d=1/4$ (blue), $d=1/16$ (green), $d=1/32$ (orange), and $d=1/128$ (purple).
    Completing all left (or right) steps first leads to catastrophic errors as the rank of the approximation increases in the experiments where the spectra of $A$ and $B$ are closest.
    Alternating between left and right steps appears to be stable in all cases. }
    \label{fig:shift_ordering}
\end{figure}

Suppose that $X_0=0$ and $C=UV^*$ with $U\in \C^{m \times \rho}$ and $V\in \C^{n \times \rho}$.
By \eqref{eq:rank-update}, we have $X_{k} = \sum_{1 \le j \le k}Z_jY_j^*$,  
    with each $Z_j \in \C^{m \times \rho}$ and $ ~ Y_j \in \C^{n \times \rho}$.
To compute these factors efficiently, 
we rewrite the Sylvester matrix equation \eqref{eq:sylv} as $\mathcal{L}X=C$, where
$\mathcal{L}(X) = AX-XB$,
and introduce the operator split, $\mathcal{L}=\mathcal{M}_k - \mathcal{N}_k$, given by
\[
    \mathcal{M}_k(X) = 
    \begin{cases}
        (A-\beta_kI)X, & \text{Left step at iteration } k,\\
        -X(B-\alpha_kI), & \text{Right step at iteration } k,
    \end{cases}
\]
and
\[
    \mathcal{N}_k(X) =
    \begin{cases}
        X(B-\beta_kI), & \text{Left step at iteration } k,\\
        -(A-\alpha_kI)X, & \text{Right step at iteration } k.
    \end{cases}
\]
We have abused the notation of the shift parameter indexing; $\alpha_k$ is not the $k$th right shift but rather the right shift at the $k$th step.
Generalized ADI is then specified by
\[
    \mathcal{M}_{k+1}(X_{k+1}) = \mathcal{N}_{k+1}(X_k) + C.
\]
Generalized ADI is a fully commuting splitting scheme, i.e., $\mathcal{M}_k$ commutes with $\mathcal{N}_j$ and $\mathcal{M}_j$ for all $(j,k)$.
Then, it follows from Proposition 2.1 of~\cite{schulze2024unifying} that 
\begin{equation}\label{eq:split_recurrence}
    X_{k+1} = X_k + \Delta X_{k+1} = X_k + \mathcal{M}_{k+1}^{-1}(R_{k}),    
\end{equation}
where the residuals, $R_j:=C-(AX_j-X_jB)$, satisfy the recurrence
\begin{equation}\label{eq:residual_recurrence}
    R_{j+1} = \mathcal{M}_{j+1}^{-1}(\mathcal{N}_{j+1}(R_{j}))=\mathcal{N}_{j+1}(\Delta X_{j+1}).
\end{equation}
Together, \eqref{eq:split_recurrence} and \eqref{eq:residual_recurrence} yield an efficient scheme for updating the iterates and residuals. 
Letting $R_j=\tilde Z_j \tilde Y_j^*$ and $\Delta X_{j}=Z_jY_j^*$ be low-rank factorizations, a left step gives
\[
Z_{j+1}Y_{j+1}^*= (A-\beta_{j+1}I)^{-1}\tilde Z_j \tilde Y_j^* ~~\text{and}~~ \tilde Z_{j+1}\tilde Y_{j+1}^* = Z_{j+1}Y_{j+1}^*(B-\beta_{j+1}I),
\]
and a right step gives
\[
Z_{j+1}Y_{j+1}^*= -\tilde Z_j \tilde Y_j^*(B-\alpha_{j+1}I)^{-1} ~~\text{and}~~ \tilde Z_{j+1}\tilde Y_{j+1}^* = -(A-\alpha_{j+1}I)Z_{j+1}Y_{j+1}^*.
\]
Grouping terms, we obtain the update rules
\begin{align*}
    &Z_{j+1} = t_{j+1}(A) \tilde Z_{j}, ~~\tilde Z_{j+1} = l_{j+1}(A)Z_{j+1},\\
    &Y_{j+1} = l^{-*}_{j+1}(B) \tilde Y_{j}, ~~\tilde Y_{j+1} = t^{-*}_{j+1}(B)Y_{j+1},
\end{align*} 
where
\begin{subequations}
\vspace{-10pt}
\begin{align}
    t_{k}(z) &= 
    \begin{cases}
        (z-\beta_{k})^{-1}, & \text{Left step at iteration } k,\\
        -1, & \text{Right step at iteration } k,     
    \end{cases}\label{eq:tk}
    \\
    l_{k}(z) &= 
    \begin{cases}
        1, & \text{Left step at iteration } k,\\
        (z-\alpha_{k}), & \text{Right step at iteration } k.
    \end{cases}\label{eq:lk}
\end{align}
\end{subequations}
We have derived the g-fADI iteration, but there is one issue to address.
If we run the iteration above with alternating left and right steps, i.e., in the setting of standard ADI,
we obtain $X_{2k} =\sum_{1 \le j \le 2k} Z_jY_j^{*}$.
Equation \eqref{eq:rank-update} indicates that $\rank(X_{2k})=k\rho$, but we are representing the final iterate with a sum of $2k$ many rank-$\rho$ outer products.
This redundancy can be rectified by explicitly compressing the recurrence relation for alternating steps. 
The compressed recurrences are given by
\begin{align}\label{eq:alternating-Z}
    Z_{k+1} &= \gamma_{k+1}(A-\beta_{k+1}I)^{-1}\tilde Z_k, ~~ \tilde{Z}_{k+1} = \frac{1}{\gamma_{k+1}}(A-\alpha_{k+1}I)Z_{k+1},
\end{align}
and
\begin{align}\label{eq:alternating-Y}
    Y_{k+1} &= (B-\alpha_{k+1}I)^{-*}\tilde Y_k, ~~
    \tilde{Y}_{k+1} = (B-\beta_{k+1}I)^*Y_{k+1},
\end{align}
where $\gamma_{k+1} = (\beta_{k+1}-\alpha_{k+1})$.
The derived iteration can proceed with arbitrary shift order, and Proposition 2.4 of~\cite{schulze2024unifying} suggests that the final result is invariant to the ordering.
This is only true in exact arithmetic, and the output can be highly sensitive to the shift ordering.
\Cref{fig:shift_ordering} shows four experiments, two where the shift ordering is inconsequential and two where it is critical.

Assuming that all the shifts need to be used to achieve the desired accuracy,
    we recommend that alternating steps are taken until one of the shift parameter lists is exhausted.
This ensures there is no storage redundancy and our experiments suggest
    it is the more numerically stable choice.
Pseudocode for g-fADI with alternating steps completed first is given by \cref{alg:g-fADI}.

\begin{algorithm}[!htbp]
\caption{Generalized factored ADI (g-fADI)}
\label{alg:g-fADI}

\Proc{\textsc{g-fADI}}
{
    \KwIn{$(A,B,U,V)$, left shifts $\{\beta_j\}_{1}^q$, right shifts $\{\alpha_\ell\}_{1}^p$}
    \KwOut{$(Z,Y,\widetilde Z,\widetilde Y),$ where $X=ZY^*$ and $\smash{\widetilde Z~\widetilde Y^*} = UV^*-(AX-XB)$}
    
    \BlankLine
    \BlankLine
    
    $\widetilde Z \gets U$,\;
    $\widetilde Y \gets V$,\;
    $j \gets 1$,\;
    $\ell \gets 1$\;

    \While{$j  \le  q$ or $\ell \le p$}{
        \uIf{$j  \le  q$ and $\ell \le p$}{
            $(\hat Z,\hat Y,\widetilde Z,\widetilde Y)
            \gets
            \AlternatingStep(
                \widetilde Z,\widetilde Y,
                \beta_{j},\alpha_{\ell},A,B)$\;

            $j \gets j+1$,\;
            $\ell \gets \ell+1$\;
        }
        \uElseIf{$j \le q$}{
            $(\hat Z,\hat Y,\widetilde Z,\widetilde Y)
            \gets
            \LeftStep(\widetilde Z,\widetilde Y,\beta_{j},A,B)$\;

            $j \gets j+1$\;
        }
        \Else{
            $(\hat Z,\hat Y,\widetilde Z,\widetilde Y)
            \gets
            \RightStep(\widetilde Z,\widetilde Y,\alpha_{\ell},A,B)$\;

            $\ell \gets \ell+1$\;
        }

        \BlankLine
            
        $Z \gets [Z, \hat{ Z}]$,\;
        $Y \gets [Y, \hat{ Y}]$
    }
    \Return $(Z,Y,\widetilde Z,\widetilde Y)$
}

\BlankLine
\hrule
\BlankLine
\textbf{Subroutines}\;

\BlankLine

\Fn{\LeftStep{$\widetilde Z,\widetilde Y,\beta,A,B$}}{
    $Z \gets (A-\beta I)^{-1}\widetilde Z$,\;
    $Y \gets \widetilde Y$\;

    $\widetilde Z \gets Z$,\;
    $\widetilde Y \gets (B^*-\bar\beta I)Y$\;

    \KwRet{$(Z,Y,\widetilde Z,\widetilde Y)$}\;
}

\BlankLine

\Fn{\RightStep{$\widetilde Z,\widetilde Y,\alpha,A,B$}}{
    $Z \gets -\widetilde Z$,\;
    $Y \gets (B^*-\bar\alpha I)^{-1}\widetilde Y$\;

    $\widetilde Z \gets -(A-\alpha I)Z$,\;
    $\widetilde Y \gets Y$\;

    \KwRet{$(Z,Y,\widetilde Z,\widetilde Y)$}\;
}

\BlankLine

\Fn{\AlternatingStep{$\widetilde Z,\widetilde Y,\beta,\alpha,A,B$}}{
    $\gamma \gets \beta-\alpha$\;

    $Z \gets \gamma(A-\beta I)^{-1}\widetilde Z$,\;
    $Y \gets (B^*-\bar\alpha I)^{-1}\widetilde Y$\;

    $\widetilde Z \gets (A-\alpha I)Z/\gamma$,\;
    $\widetilde Y \gets (B^*-\bar\beta I)Y$\;

    \KwRet{$(Z,Y,\widetilde Z,\widetilde Y)$}\;
}

\end{algorithm}

\subsection[How many left or right steps?]{How many left or right steps?\nopunct} Depending on the sizes of $A$ and $B$, as well as their spectral properties, it can be advantageous to use differing numbers of left and right steps. This is precisely when g-fADI is more effective than fADI.
Let $T_A,T_B$ be the costs of shifted inverses and $M_A,M_B$ the costs of matrix-vector products. The computational complexity of taking a left step is $O(\rho(T_A+M_B))$ and the computational complexity of taking a right step is $O(\rho(T_B+M_A))$.
Letting $C_L$ and $C_R$ denote the costs of a left and right steps, respectively, the computational cost of executing g-fADI with $q$ left-steps and $p$ right-steps is $O(q C_L + pC_R).$
Letting $N=(p+q)$ and $\omega=p/(p+q)$, we rewrite this cost as $O(N(\omega C_R + (1-\omega)C_L))$. Using \cref{thm:asymptotic} along with the error bound \eqref{eq:generalized-ADI-error-bound-nonnorm}, we have that for a given $\omega \in [0,1]$, the total number of iterations required to reach a tolerance $\varepsilon$ is
\[
N_\varepsilon \sim \frac{\log(1/\varepsilon)}{\mathscr{F}(\omega)},
\]
where $\mathscr{F}(\omega)$ is a nonnegative, concave function defined by the geometry of the spectral sets of $A$ and $B$ and introduced by Levin and Saff in~\cite{levin1994optimal}. 
The equation above does not explicitly show any dependence on the problem size $(m,n)$.
In many settings, the spectra of $A$ and $B$ slowly approach one another as $m$ and $n$ increase.
This has the effect of decreasing $\mathscr{F}(\omega)$, hence increasing the number of iterations required to reach a certain tolerance.

For a fixed ratio of left to right steps, the computational complexity of g-fADI is 
\begin{equation}\label{eq:g-fADI-cost}
O\left(\frac{\rho \log\left(\frac{1}{\varepsilon}\right)}{\mathscr{F}\left(\omega\right)}\left(\omega C_R + \left(1-\omega\right)C_L\right)\right).
\end{equation}
Classical ADI corresponds to $\omega=1/2$, and it is not hard to find examples where $\mathscr{F}(1) > \mathscr{F}(1/2)$.
In these settings, one-sided g-fADI ($q=0$) is more efficient in step count than classical fADI ($p=q$).
If also $C_R \ll C_L$, then the savings will be considerable.
While optimizing \eqref{eq:g-fADI-cost} is outside the scope of this paper, we give several examples in \cref{sec:examples} that show when and how $p \neq q$ can be advantageous.

\subsection{One-sided g-fADI}\label{subsec:onesided-gfadi}

The g-fADI algorithm can be run with only right (or only left) steps. Suppose that only right steps are taken, with shifts $\{\alpha_i\}_{i=1}^p$. Then $t_k\equiv-1$ and $l_k(z)=z-\alpha_k$ in~\eqref{eq:tk}--\eqref{eq:lk}, so that the $k$th iterate is given by
\[
X_k=-\sum_{j=1}^k\left(\prod_{\ell=1}^{j-1}(A-\alpha_\ell I)\right)UV^*\left(\prod_{\ell=1}^{j}(B-\alpha_\ell I)^{-1}\right).
\]
The iteration never solves a shifted system with $A$, since every step costs one shifted solve with $B$ and one product with $A$. 
In the extreme case when $B=\beta$ is a scalar, the one-sided iteration is the nonstationary Richardson iteration $x^{(i)}=x^{(i-1)}+(u-(A-\beta I)x^{(i-1)})/(\alpha_i-\beta)$ for $(A-\beta I)x=u$~\cite{Richardson1911}. With each $\alpha_i=0$ (far from optimal), it computes a truncated Neumann series, and with Chebyshev zeros and real $\beta$, it is the Chebyshev--Richardson iteration~\cite{Young1953}. To our knowledge, the connecgtion to ADI is new, as is the extension to nondiagonal and nonnormal $B$. %

As a simple illustration, suppose $A$ is a dense symmetric $m\times m$ matrix with spectrum in $[1/4,1]$, let $B$ be a symmetric $200\times 200$ matrix with spectrum in $[-1,-1/4]$, and let $UV^*$ be a rank 1 matrix. With $k$ shifts selected as the Chebyshev points of the first kind (the roots of $T_k$) rescaled to the interval $[1/4,1]$, the relative error is at most $1/T_k(7/3)$, where $T_k$ is the $k$th Chebyshev polynomial. 
Lemma~\ref{lem:asym_rates} shows that these shifts attain the optimal rate.

It takes g-fADI $23$ steps to reach an error of $10^{-14}$. Classical fADI with optimal shifts needs about as many shifted solves, but half of them are with $A$, so its cost grows like $m^3$. The one-sided iteration needs only products with $A$, and at $m=8000$, it is $40$ times faster than classical fADI and $125$ times faster than Bartels--Stewart.

We analyze when one-sided iterations are viable in \cref{sec:one-sided}, and we return to them in \cref{sec:examples}, where $A$ is a discretized three-dimensional operator and the shifts come from rigorous bounds on its spectrum. 

%%%%%%%%%%%%%%%%%%%%%%%%%%%%%%%%%%%%%%%%%%%%%%%%%%%%%%%%%%%%%%%%%%%%%%%%%%%%%%%%%%%%%%%%%%%%%%%%%%%%%%%%%%%%%
%------------------------COMPUTING SHIFTS-------------------------------------------------------------------%
%%%%%%%%%%%%%%%%%%%%%%%%%%%%%%%%%%%%%%%%%%%%%%%%%%%%%%%%%%%%%%%%%%%%%%%%%%%%%%%%%%%%%%%%%%%%%%%%%%%%%%%%%%%%%
\section{Convergence rates and computing shift parameters}\label{sec:theory}
The shifts enter the error bound~\eqref{eq:generalized-ADI-error-bound-nonnorm} through Zolotarev's third problem on $E=W(A)$ and $F=W(B)$ (or $E = \Lambda(A), F = \Lambda(B)$ when $A, B$ are normal), which is solved in closed form only for a handful of geometries, even when $p=q$. In this section, we review common shift-selecting practices and propose a framework for computing them in the generalized setting, motivated by the asymptotic results of~\cite{levin1994optimal}.

The shifts do not require $W(A)$ and $W(B)$ themselves, only sets that contain them. If $B$ is small, its field of values can be estimated or computed directly~\cite{Johnson1978}. For a discretized differential operator, the spectrum or the field of values can often be enclosed by analysis of the operator itself or by inexpensive rigorous bounds, as in the examples of \cref{sec:examples}. When no useful enclosure is available, projection-based methods are usually the better choice. Our purpose here is to understand the convergence of g-ADI when such sets are available.

\subsection{Shift computation methods} Several heuristics for computing suboptimal shift parameters have been developed in the literature.
One option is to approximate $E$ and $F$ by larger sets for which the problem is solved, often disks or intervals, and use the closed-form shifts of the enclosing configuration~\cite{EllnerWachspress1991,LuWachspress1991,Wachspress1990}.
The convergence rate can degrade severely with a poor enclosure.
Alternatively, one can use Leja--Bagby points~\cite{Bagby1967, Leja1957}, which are generated by a greedy recursion directly from $E$ and $F$: having chosen the first $k$ zeros and poles, the next zero is placed where the current rational function is largest on $E$, and the next pole where it is smallest on $F$.
Each step is a scalar maximization over a discretization, so the construction is cheap and applies to arbitrary compact sets; this is the form in which Leja points enter ADI shift selection~\cite{BennerKurschnerSaak2014etna, LevenbergReichel1993}.
The resulting point distributions converge to the equilibrium measures of the condenser $(E,F)$~\cite{Bagby1967}, recovering the asymptotically optimal rate of Gonchar~\cite{Gonchar1969}.
This is a statement about the limit, however, and ADI is run at finite $p$ and $q$.
The gap between the Leja--Bagby rate and the Zolotarev rate at finite $p,q$ can be substantial, and is not controlled by the asymptotic theory.

Most recently, there have been developments in data-driven methods for Zolotarev problems, built on the AAA--Lawson algorithm for minimax rational approximation~\cite{NakatsukasaTrefethen2020} and independently on the Loewner framework~\cite{poussotvassal2025zolotarev}.
These methods solve Zolotarev's third problem by using its equivalence with Zolotarev's fourth
problem~\cite{istace1995third}, a sign-function rational approximation problem.
These methods are able to compute the optimal shifts on various difficult geometries.
Nearness to optimality can be verified by checking that the error equioscillates in the real case and is near-circular on $\partial E$ and $\partial F$, in the complex setting~\cite{achieser2013theory,Starke1992}.
In this section, we discuss the convergence properties of generalized Zolotarev numbers, and we generalize the AAA-based Zolotarev algorithm of~\cite{trefethen2024computation} for the case where $p \neq q$.

%%%%%%%%%%%%%%%%%%%%%%%%%%%%%%%%%%%%%%%%%%%%%%%%%%%%%%%%
%%%%%%%%%%%%%%%%%%%%%%%%%%%%%%%%%%%%%%%%%%%%%%%%%%%%%%%%
\subsection{Approximation results}
Our primary motivation for considering the type-$(p,q)$ Zolotarev ratio problem is to minimize the cost of g-fADI, given by~\eqref{eq:g-fADI-cost}.
The two crucial factors in the complexity are the costs of left/right steps and the total number of iterations required. 
Our goal is to choose $p$ and $q$ to optimize the total number of flops required to achieve a specified accuracy. 
In this section, we review the theory in~\cite{levin1994optimal}, which characterizes the decay of type-$(p,q)$ Zolotarev numbers, then use it to understand how the convergence of g-fADI depends on the ratio of left steps to right steps.

Even when the costs of left and right steps are the same, it can be technically more \textit{efficient} to have $p \ne q$ because of the decay behavior of the generalized Zolotarev numbers.
In the context of g-fADI, the implication is that the target accuracy can be achieved in fewer iterations using a disproportionate amount of left and right steps.
Levin and Saff explore type-$(p,q)$ Zolotarev numbers at length in~\cite{levin1994optimal}, and
    their asymptotic results (see \cref{thm:asymptotic}) 
    inform and motivate the present work.
    
\begin{figure}[t!]
    \centering
    \includegraphics{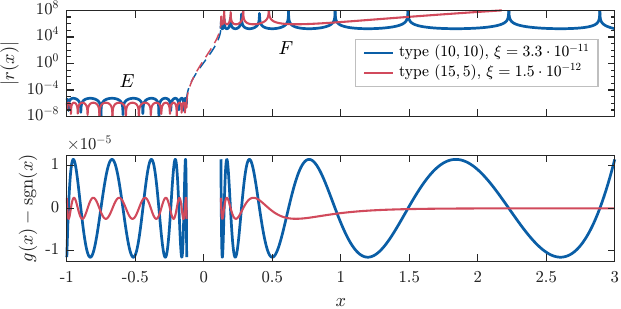}
    \caption{Zolotarev rationals for $E=[-1,-1/8]$ and $F=[1/8,3]$ of types $(10,10)$ (blue) and $(15,5)$ (red) are shown in the top plot.
    We have $\xi_{10,10}\approx \num{3.3e-11}$ and $\xi_{15,5}\approx \num{1.5e-12}$.
    The Zolotarev rationals were computed via \cref{alg:modified aaa} in \cref{app:z3alg} with Lawson iteration applied afterwards.
    The corresponding sign approximation error is shown in the bottom plot, with the equioscillation indicating optimality of the approximants in each case.}
    \label{fig:optimal_ray_motivation}
\end{figure}

\Cref{fig:optimal_ray_motivation} shows a concrete example where the optimal type-$(15,5)$ rational function achieves a modestly smaller Zolotarev ratio than the optimal type-$(10,10)$ rational function.
The total number of poles and zeros is 20 in each case and $\xi_{10,10} \approx 20\xi_{15,5}$.
For a Sylvester matrix equation with $\Lambda(A)$ and $\Lambda(B)$ matching the geometries in the experiment below, g-fADI (using 15 right and 5 left steps) would have an error bound $20$ times smaller than that of fADI (using 10 right and 10 left steps).
It is entirely possible that there is an even more efficient choice of $(p,q)$ for this geometry.

We have selected an example with only a modest improvement                   
    in the lopsided case for legibility of the plots.
The improvement can be larger for other geometries and choices of $(p,q)$. 
To understand how much larger, we turn to the work of Levin and Saff \cite{levin1994optimal}, 
    where they generalize various properties of classical Zolotarev numbers.
A key property of classical Zolotarev numbers is that they decay exponentially with rate depending on $(E,F)$, i.e.,
\[
    \lim_{q \to \infty}(\xi_{q,q})^{1/(2q)} = \exp\bigg(\frac{-1}{2\text{cap}(E,F)}\bigg),
\]
where $\text{cap}(E,F)$ is the capacitance of the condenser $(E,F)$~\cite{Gonchar1969} .
Particularly, the distribution of poles and zeros of the Zolotarev rationals
    as $k\to \infty$ is described by the signed equilibrium measure on the condenser $(E,F)$.
Levin and Saff extend this analysis to the type-$(p,q)$ case by considering
    signed equilibrium measures disproportionately concentrated on $E$ and $F$.
The asymptotic distribution of poles and zeros of a type-$(p,q)$ Zolotarev rational
    is described by a signed equilibrium measure $\mu$ on $(E,F)$ having $|\mu(E)|/|\mu(F)|= p/q$~\cite{levin1994optimal}.

This technique is used to show that type-$(p,q)$ Zolotarev numbers also decay exponentially, 
    now with rate depending on $E,F,$ and $\lambda=p/q$.
The behavior of type-$(p,q)$ Zolotarev numbers is made precise in \cref{thm:asymptotic} below,
    which paraphrases the results in \cite{levin1994optimal} most relevant to the present work.
\begin{theorem}[Levin and Saff~\cite{levin1994optimal}]\label{thm:asymptotic}
    Let $E,F \subset \C$ be disjoint, compact sets each having positive logarithmic capacity,
    and let $\{(m_k,n_k)\}_{k=1}^\infty$ be a sequence of ordered pairs of nonnegative integers
    satisfying
    \begin{equation}\label{eq:ray}
        \lim_{k \to \infty}\frac{m_k}{n_k} = \lambda, ~~
        \lim_{k \to \infty} m_k +n_k=  \infty,
    \end{equation}
    where $0 < \lambda < \infty$.
    Then,
    \begin{equation}\label{eq:saff-decay}
        \lim_{k \to \infty} \xi_{m_k,n_k}^{1/(m_k+n_k)} = e^{-\mathscr{F}(\omega)},
    \end{equation}
    where $\omega=\lambda/(1+\lambda)$ and $\mathscr{F}(\cdot)$ is a nonnegative, concave function.
\end{theorem}
\begin{proof}
    \Cref{thm:asymptotic} combines Theorems 3.2 and 6.1 of~\cite{levin1994optimal},
    the proofs of which can be found there.
\end{proof}
\begin{remark}
    The function $\mathscr{F}(\omega)$ is determined by $(E,F)$ and 
        is given by~\cite[eq.~(3.6)]{levin1994optimal}.
    We omit the definition here since our discussion requires only the property
        that $\mathscr{F}(\omega)$ is nonnegative and concave.
\end{remark}

The nonnegativity and concavity of $\mathscr{F}(\omega)$ in \cref{thm:asymptotic} 
    implies a limit on how much can be gained asymptotically from lopsided rationals, particularly,
    \begin{equation}\label{eq:max_benefit}
        \frac{\mathscr{F}(\omega)}{\mathscr{F}(1/2)} \le 1+2|\omega-1/2|, \;\; \forall \omega \in (0,1).
    \end{equation}
    An immediate consequence
    is that the maximum asymptotic benefit from $p \ne q$ is a factor of 2 in the exponential decay rate\footnote{We thank Nick Trefethen and Yuji Nakatsukasa for positing this question to us at the BIRS workshop on computational rational approximation.}.
    Example 5.1 of~\cite{levin1994optimal} shows that in the case of $E$ and $F$ being concentric circles, the factor of 2 is attained.
    The implication is that, for some geometries, g-fADI requires half as many iterations as fADI, with the maximal savings achieved by a one-sided iteration.

    While halving the number of required iterations can be a significant improvement in some settings, the biggest savings in g-fADI come from avoiding expensive iterations altogether.
    Many applications have that $C_L$ and $C_R$ differ in magnitude, often due to one matrix having significantly cheaper shift-inverts.
    In these scenarios, the focus is not for lopsidedness to decrease the iteration count but rather to greatly cheapen the per-iteration cost without increasing the iteration count by too much.
    For example, if right steps are much cheaper, i.e., $C_R \ll C_L$, one-sided iterations will be significantly faster so long as
    \[
        \frac{C_R}{C_L} \ll \frac{\mathscr{F}(1)}{\mathscr{F}(1/2)}.
    \]
    On some geometries, this condition is guaranteed to fail.
    In Example 5.2 of~\cite{levin1994optimal}, Levin and Saff analyze the geometry $E=\{|z|=r<1\}\cup\{|z|=R>1\}$ and $F=\{|z|=1\}$, finding $\mathscr{F}(0)=\mathscr{F}(1)=0$ and $\mathscr{F}(1/2) \ne 0$.
    In this case, one-sided iterations diverge and so offer no benefit, regardless of $C_R$ and $C_L$.
    While $p \ne q$ can improve the decay rate by up to a factor of two, it can be arbitrarily worse than $p = q$.
    For this reason, it is important to assess whether a given geometry is amenable to $p \ne q$. We now turn to the development of a computational method for doing that. 
    
\subsection{Constructing generalized Zolotarev rationals via an equivalence theorem}\label{sec:computing shifts}
 Zolotarev's third problem is equivalent to a minimax rational approximation problem, called Zolotarev's fourth problem, in the classical setting and under mild conditions on $E$ and $F$~\cite{istace1995third}.
In the generalized setting, Z3 is equivalent to a closely related approximation problem.

\smallskip
\begin{mdframed}[nobreak=true]
\noindent \textbf{(Z4) Zolotarev's fourth problem.} Given disjoint compact sets $E$ and $F$ in the complex plane and fixed nonnegative integers $p,q$, find a rational function $g \in \mathcal{R}_{p,q}$ that minimizes
\begin{equation}\label{eq:Z4_error}
    \tau(g) := \max_{z \in E \cup F}\bigl|g(z) - \sgn_{E/F}(z)\bigr|,
    \qquad
    \sgn_{E/F}(z) :=
    \begin{cases}
        -1, & z \in E,\\
        +1, & z \in F.
    \end{cases}
\end{equation}
The minimum is denoted by
\begin{equation}\label{eq:Z4_number}
    \tau_{p,q}(E,F) := \min_{g \in \mathcal{R}_{p,q}} \tau(g).
\end{equation}
\end{mdframed}
\smallskip

\Cref{thm:lopsided-equiv} gives an abridged version of the generalized equivalence, providing the transformations that map solutions of one problem to the other.
In light of the equivalence, any minimax rational approximation algorithm can be used to solve Z3 in the classical setting.
In~\cite{trefethen2024computation}, AAA--Lawson was used to successfully compute Zolotarev rationals on a variety of difficult geometries.
Solving Z3 in the generalized setting is not as straightforward, due to a certain constraint in the equivalent problem.
Consequently, it is not feasible to use, without modification, preexisting minimax rational approximation algorithms to solve generalized Z3.

\begin{theorem}[generalized Zolotarev equivalence]\label{thm:lopsided-equiv}
Let $E,F \subset \C$ be disjoint, compact sets and $p,q \in \mathbb{N}$.
Let  $k=\max\{p,q\}$.
Define the transform
    \begin{equation}\label{eq:S}
        S[g](z) = \frac{(1-\xi)/(1+\xi)+ g(z)}{(1-\xi)/(1+\xi) -g(z)},
        ~~
        \xi = \Bigg(\frac{\tau(g)}{1+\sqrt{1-\tau(g)^2}}\Bigg)^2,
    \end{equation}
     where $\tau(g)$ is as in~\eqref{eq:Z4_error}.
    Define the set
    \begin{equation}\label{eq:G}
        \G = \{g \in \RR_{k,k} ~|~ S[g] \in \RR_{p,q}\}.
    \end{equation}
    Suppose $g^*$ solves the optimization problem
    \begin{equation}\label{eq:gstar}
        g^* = \underset{g \in \G}{\text{argmin}}~\tau(g).
    \end{equation}
    Then, $r^*=S[g^*] \in \RR_{p,q}$ solves
    the type-$(p,q)$ third Zolotarev problem.
    Particularly, 
    \begin{equation}\label{eq:rstar}
    \xi(r^*) = \xi_{p,q} = \Big(\frac{\tau(g^*)}{1+\sqrt{1-\tau(g^*)^2}}\Big)^2.
    \end{equation}
\end{theorem}
\begin{proof}
    A proof can be found in~\cite{LeBaillyThiran2000}.  
\end{proof}
\noindent 
In the setting where $p=q$, we have that $\G=\RR_{p,p}$ in \eqref{eq:G}.
This simplifies the equivalence greatly: problem \eqref{eq:gstar} is a standard type-$(p,p)$ minimax approximation problem.
In the generic case, solving \eqref{eq:gstar} exactly is much harder since one must optimize over $\G$ given by \eqref{eq:G}, a set of functions defined by a nonlocal condition. 
We address this difficulty by replacing $S$ with the more tractable transformation, 
    \begin{equation} 
    \label{eq:hats}
    \tilde{S}[g](z) = \frac{1 + g(z)}{1-g(z)},
    \end{equation}
    which is a natural approximation to consider since we are most interested in settings where $0 < \xi \ll 1$.   
In \cref{lem:approx_S}, we characterize how close $\tilde{S}[g](z)$ is to $S[g](z)$.
This allows us to interpret $\tilde{S}[\cdot]$ as defining an approximate equivalence and use it to build an algorithm for computing generalized Zolotarev rationals.
\begin{lemma}\label{lem:approx_S}
Suppose $E,F \subset \C$ are disjoint and compact, 
    that $p,q \in \mathbb{N}$, and let $\tilde S[\cdot]$
    be as in \eqref{eq:hats}.
Denote $k = \max\{p,q\}$ and let 
\begin{equation}\label{eq:gtilde_set}
    \Gt = \{g \in \RR_{k,k} ~|~ \tilde S[g] \in \RR_{p,q}\} \subset \RR_{k,k}.
\end{equation}
If $\tilde{g} \in \Gt$ solves the optimization problem
\begin{equation}\label{eq:gtilde_problem}
    \tilde{g} = \text{argmin}_{g \in \Gt}\tau(g),
\end{equation}
then the function $\tilde{r}(z)=\tilde S[\tilde{g}](z) \in \RR_{p,q}$ 
has Zolotarev ratio bounded by
\begin{equation}
    \xi(\tilde r) 
    \le \xi_{p,q} + O(\xi_{p,q}^{3/2}), ~~\text{for } \xi_{p,q} \ll 1,
\end{equation}
where $\xi_{p,q}$ is the
    type-$(p,q)$ Zolotarev number.
\end{lemma}
\begin{proof}
    Let $g^*\in\G$ solve~\eqref{eq:gstar}, and set $\hat g=g^*/\gamma$ with $\gamma=(1-\xi_{p,q})/(1+\xi_{p,q})$. By~\eqref{eq:G},
    \[
        \tilde S[\hat g]=\frac{1+g^*/\gamma}{1-g^*/\gamma}=S[g^*]\in\RR_{p,q},
    \]
    so $\hat g\in\Gt$. Since $\max_{E\cup F}|g^*|\le1+\tau(g^*)$ and $1/\gamma-1=2\xi_{p,q}/(1-\xi_{p,q})$, the solution $\tilde g$ of~\eqref{eq:gtilde_problem} satisfies
    \[
        \tau(\tilde g)\le\tau(\hat g)\le\tau(g^*)+\frac{2\xi_{p,q}}{1-\xi_{p,q}}\bigl(1+\tau(g^*)\bigr)=\tau(g^*)+O(\xi_{p,q}).
    \]
    Since $\tau\mapsto(\tau/(2-\tau))^2$ is smooth and increasing on $[0,1]$, we have
    \begin{align*}
        \xi\bigl(\tilde S[\tilde g]\bigr)
        &\le\Bigl(\frac{\tau(\tilde g)}{2-\tau(\tilde g)}\Bigr)^2
        \le\Bigl(\frac{\tau(g^*)}{2-\tau(g^*)}\Bigr)^2+O\bigl(\tau(g^*)\xi_{p,q}\bigr)\\
        &=\Bigl(\frac{\tau(g^*)}{1+\sqrt{1-\tau(g^*)^2}}\Bigr)^2+O\bigl(\tau(g^*)^3\bigr)+O\bigl(\tau(g^*)\xi_{p,q}\bigr)
        =\xi_{p,q}+O(\xi_{p,q}^{3/2}),
    \end{align*}
    where the last step uses~\eqref{eq:rstar}, which also gives $\tau(g^*)=O(\xi_{p,q}^{1/2})$.
\end{proof}

We have established that $\tilde S[\cdot]$ transforms solutions of~\eqref{eq:gtilde_problem} to near-solutions of the generalized Zolotarev ratio problem.
The AAA--Lawson framework as used in~\cite{trefethen2024computation} is particularly amenable to searching over $\Gt$ as given by~\eqref{eq:gtilde_set}.
The AAA--Lawson algorithm computes minimax rational approximants in barycentric form, 
\[
g(z) = \frac{\sum_{j=1}^{k+1} a_j/(z-s_j)}{\sum_{j=1}^{k+1} b_j/(z-s_j)}.
\]
A simple calculation establishes that
\[
\tilde S[g](z) = \frac{\sum_{j=1}^{k+1}(b_j + a_j)/(z-s_j)}{\sum_{j=1}^{k+1}(b_j - a_j)/(z-s_j)},
\]
showing that $\tilde S[\cdot]$ is linear in the barycentric weights.
Enforcing $\tilde S[g]\in \mathcal{R}_{p,q}$ amounts to constraining the barycentric weights to a subspace.
Pseudocode for the AAA portion of the described framework, including the constraint on the barycentric weights, is given in \cref{alg:modified aaa} in \cref{app:z3alg}.
The output of \cref{alg:modified aaa} can then be used as the initialization for a Lawson iteration in order to obtain the minimax solution.
\Cref{fig:zolo_gallery} shows three experiments with Zolotarev rational functions of different types on various geometries.
Using this computational method is cheap and allows us to explore and understand the behavior of $\xi_{p,q}$ on a variety of geometries. It also allows one to cheaply compute near-optimal g-ADI shift parameters (the poles and zeros of the computed Zolotarev rational of type $(p,q)$).%

One may notice, however, that in the case where $p$ and $q$ diverge substantially, it can become numerically unstable to compute the constraint conditions on the barycentric weights. This is unfortunate since one of the most useful versions of g-fADI is the case where $p=0$ or $q=0$. We supply some analysis for common geometries (e.g., intervals) in the next section. In the examples of \cref{sec:examples}, we therefore do not compute the shifts with AAA--Lawson. Instead, we enclose the relevant spectrum in an interval whose endpoints come from inexpensive but rigorous eigenvalue bounds, take the shifts to be Chebyshev points on that interval, and choose the number of steps from the resulting explicit error bound. When the matrices are far from normal, we enclose their numerical ranges instead and pay a factor of at most $4$~\cite{jin2026numerical,lorist2026}. The cost of computing these bounds is included in every reported timing.

\subsection{One-sided iterations}\label{sec:one-sided}
The g-fADI algorithm can offer the most speedup when the costs of left and right steps differ greatly.
In this setting, massive savings are realized by computing a one-sided iteration that avoids expensive inverts altogether.
The crucial assumption is that the number of iterations does not grow too drastically due to the completely unbalanced iteration.
Whether this assumption is reasonable or not depends on the geometries of $E$ and $F$.
In this section, we analyze the viability of one-sided iterations on problems with geometries given by asymmetric real intervals.

Real intervals are historically significant to Zolotarev's problems.
Yegor Ivanovich Zolotarev posed and solved his third and fourth problems on \textit{symmetric} real intervals of the form $E=[-1,-\delta]$ and $F=[\delta,1]$.
Since Zolotarev's problems are invariant under dilation and shifting transformations, this form represents all cases of finite-length, disjoint same-length real intervals. A natural extension to study in the type-$(p,q)$ setting is the case of asymmetric real intervals. Appealing to intuition, it might seem \textit{more efficient} to have more zeros than poles if $E$ is much smaller than $F$.

Levin and Saff performed a rigorous analysis of the generalized Zolotarev problem in the case of asymmetric intervals. In this section, we analyze the efficiency of polynomials and their reciprocals in comparison to classical Zolotarev rationals. As long as there is not too drastic a gap in efficiency, one-sided g-fADI is viable. We will work exclusively with polynomial approximants, corresponding to right-step-only g-fADI iterations. Left-step-only iterations correspond to inverse polynomial approximants, and simply to swapping $E$ and $F$ in Zolotarev's problem.

Since $\xi_{p,q}$ is invariant under dilations and shifts, we assume without loss of generality that $E=[-1,1]$ and $F=[1+d,1+d+L]$ with $d,L>0$. This form emphasizes that the important quantities are the separation $E$ and $F$ and their relative lengths. In the notation of \cref{thm:asymptotic}, \cref{lem:asym_rates} gives the decay rates of polynomials,  corresponding to $\mathscr{F}(1)$, and of classical Zolotarev rationals, corresponding to $\mathscr{F}(1/2)$, on these sets. To be clear, these facts follow straightforwardly from known results on Chebyshev polynomials and type $(n,n)$ Zolotarev rationals.

\begin{lemma}\label{lem:asym_rates}
Let $E=[-1,1]$ and $F=[1+d,1+d+L]$ with $d,L>0$. Then
\[
    \lim_{m\to\infty}\xi_{m,0}^{1/m} = e^{-\mathscr{F}(1)}
    \quad\text{and}\quad
    \lim_{n\to\infty}\xi_{n,n}^{1/(2n)} = e^{-\mathscr{F}(1/2)},
\]
where
\begin{equation}\label{eq:cheby_vs_kk}
    \mathscr{F}(1) = \log\Big(1+d+\sqrt{d^2+2d}\,\Big)
    \quad\text{and}\quad
    \mathscr{F}(1/2) = \frac{\pi K(k')}{2K(k)}.
\end{equation}
Here $k^2 = 2L/\big((2+d)(d+L)\big)$ is a cross-ratio of the endpoints of $E$ and $F$, $k'=\sqrt{1-k^2}$, and $K(k) = \int_0^1 \big((1-t^2)(1-k^2t^2)\big)^{-1/2}\,\text{d}t$ is the complete elliptic integral of the first kind.
\end{lemma}
\begin{proof}
Set $\rho = 1+d+\sqrt{d^2+2d} > 1$, so that $\rho^{-1} = 1+d-\sqrt{d^2+2d}$.
On $[1,\infty)$, the Chebyshev polynomial $T_m(x) = \tfrac12\big((x+\sqrt{x^2-1})^m+(x-\sqrt{x^2-1})^m\big)$ is increasing, and $T_m(1+d) = \tfrac12(\rho^m+\rho^{-m})$.
By the extremal growth property of Chebyshev polynomials~\cite{Rivlin1990}, every polynomial $p$ of degree at most $m$ satisfies $|p(x)| \le T_m(x)\max_{E}|p|$ for $x \ge 1$. Hence
\[
    \min_{z\in F}|p(z)| \le |p(1+d)| \le T_m(1+d)\max_{x\in E}|p(x)|,
\]
with equality throughout when $p=T_m$, since $\max_E|T_m| = 1$ and $\min_F|T_m| = T_m(1+d)$.
Therefore $\xi_{m,0} = 1/T_m(1+d) = 2/(\rho^m+\rho^{-m})$, and $\smash{\xi_{m,0}^{1/m}} \to \rho^{-1}$.

For type $(n,n)$, $\smash{\xi_{n,n}^{1/(2n)}} \to \exp\big(-1/(2\,\text{cap}(E,F))\big)$, where $\text{cap}(E,F)$ is the capacitance of the condenser $(E,F)$~\cite{Gonchar1969}.
For two disjoint real intervals, $\text{cap}(E,F) = K(k)/(\pi K(k'))$~\cite{BeckermannTownsend2019}, which gives $\mathscr{F}(1/2)$.
\end{proof}

\begin{figure}[tbp]
        \centering
        \sbox0{\scriptsize$\mathscr{F}(1)/\mathscr{F}(1/2)$}%
        \sbox1{\includegraphics{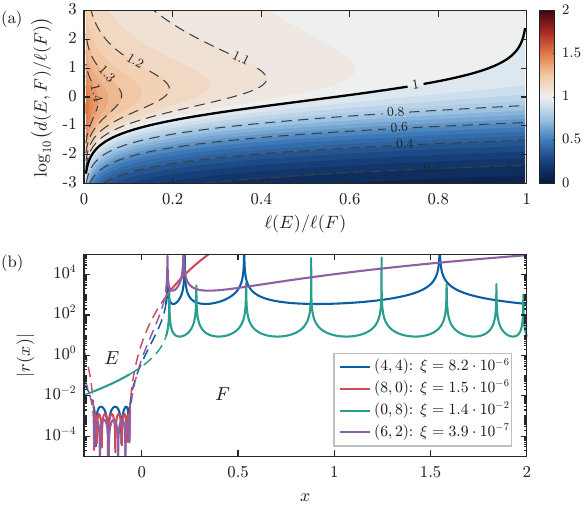}}%
        \usebox1\hspace{2pt}\raisebox{\dimexpr0.81\ht1-0.5\wd0\relax}{\rotatebox{90}{\usebox0}}
        \caption{(a) Plot of the relative efficiency of type-$(p,0)$ versus type-$(p,p)$ solutions of the Zolotarev ratio problem.
        The quantities $\ell(E)$ and $\ell(F)$ are the lengths of $E$ and $F$, respectively, and $d(E,F)$ is the separation of the two sets.
        The solid black line shows when both options are equally efficient, and the dashed black lines show various level sets of the efficiency ratio, $\mathscr{F}(1)/\mathscr{F}(1/2)$.
        A ratio of 0.2 indicates that a one-sided iteration (right steps only) would require 5 times as many shift inverts as a balanced iteration to reach the same error level.
        If the relative cost of right steps is much cheaper than that of left steps, the one-sided iteration will be much faster.
        Ratios greater than 1.0 indicate that the one-sided iteration requires \textit{fewer} total shift inverts to reach the same error level.
        (b) Zolotarev rationals for $E=[-1/4,-1/16]$ and $F=[1/8,2]$ of types $(4,4)$ (blue), $(8,0)$ (red), $(0,8)$ (teal), and $(6,2)$ (purple).
        We have $\xi_{0,8}\approx \num{1.4e-2}$, $\xi_{4,4}\approx \num{8.2e-6}$, $\xi_{8,0}\approx \num{1.5e-6}$,
        and $\xi_{6,2}\approx \num{3.9e-7}$.}
        \label{fig:asymm_intervals_analysis}
\end{figure}

% ----------------------------------------------------------------------------------------------
% ----------------------------------------------------------------------------------------------
Whenever $\mathscr{F}(1) > \mathscr{F}(1/2)$, Chebyshev polynomials are \textit{more} efficient at minimizing the Zolotarev ratio than type-$(p,p)$ rational functions.
A plot comparing the expressions in \eqref{eq:cheby_vs_kk} over various relative interval lengths 
    and separation distances is shown in \cref{fig:asymm_intervals_analysis}.

We are not claiming polynomials are the \textit{most} efficient choice of approximant.
In fact, Levin and Saff show\footnote{They do not obtain a formula for $\omega^*$.} that the optimal asymptotic ratio $\omega^* = p/(p+q)$ in this setting has $1/2<\omega^*<1$ \cite{levin1994optimal}.
\Cref{fig:asymm_intervals_analysis} shows an example where polynomials are suboptimal but still offer significant improvement over the classical Zolotarev rational. In this experiment, type-$(6,2)$ was most efficient, followed by type-$(8,0)$, both outperforming the classical type-$(4,4)$ Zolotarev rational. The type-$(0,8)$ rational achieved the worst Zolotarev ratio, and at first glance, this might seem to indicate that left-step-only iterations are a bad choice for this geometry. However, this is misleading. Approximating the decay rates,
$\mathscr{F}\big(p/(p+q)\big) \approx -\ln(\xi_{p,q})/(p+q),$
we obtain
$\mathscr{F}(0) \approx 0.54,
\mathscr{F}(1/2) \approx 1.46,
\mathscr{F}(3/4) \approx 1.84,~\text{and}~
\mathscr{F}(1) \approx 1.68.$
This shows that using left steps only will require no greater than 3 times as many iterations as alternating between left and right steps.
If left steps are significantly cheaper than right steps, this increase in iteration count is not worrisome.

%%%%%%%%%%%%%%%%%%%%%%%%%%%%%%%%%%%%%%%%%%%%%%%%%%%%%%%%%%%%%%%%%%%%%%%%%%%%%%%%%%%%%%%%%%%%%%%%%%%%%%%%%%%%%
%------------------------EXAMPLES AND NUMERICS--------------------------------------------------------------%
%%%%%%%%%%%%%%%%%%%%%%%%%%%%%%%%%%%%%%%%%%%%%%%%%%%%%%%%%%%%%%%%%%%%%%%%%%%%%%%%%%%%%%%%%%%%%%%%%%%%%%%%%%%%%
\section{Applications and examples} \label{sec:examples}

In this section, we supply several examples showcasing applications where g-fADI is useful. In each of them, inverts with $A$ are undesirable for various reasons, whereas $B$ is small or structured in a way that makes inverts very cheap (e.g., diagonal or banded).  Example 2 (\cref{sec:ex:pde}) introduces a specialized compression scheme that exploits the difference in size between the low-rank factors $Z$ and $Y$ of $X\approx ZY^*$ in memory-intensive settings. Throughout, we also supply analysis as needed to improve shift selection, modify discretizations, and better estimate error bounds when $A$ or $B$ is nonnormal with a potentially large field of values.  Details on the matrices we use can be found in \cref{app:matrices}. 

In all of the examples, the shifts are chosen a priori via the bound~\eqref{eq:onesided} of \cref{prop:onesided}, or its extension to a region that contains $W(A)$ (\cref{app:switch}), so that the error bound is at most $\varepsilon=10^{-8}$.  We find that the bounds are somewhat pessimistic in practice: our errors are below this threshold, and often well below it. All timings are wall-clock times in MATLAB R2025b on a laptop with an 8-core Intel Core Ultra 7 258V and 32~GB of memory.
All of our algorithms and experiments are publicly available on GitHub~\cite{gfadiCode}.

\subsection{Example 1: A basic test case with large $A$, small $B$ }\label{sec:tb:setup}
We first consider a basic test problem,  $AX-XB=UV^*$, where $A$ is the second-order finite-difference Laplacian~\eqref{eq:lap}, and $B=-L^T$, where $L\in\mathbb{R}^{200\times200}$ is the second-order finite-difference matrix~\eqref{eq:tbL} of $(1-\epsilon\,\partial_t^2)$ with $\epsilon\approx0.062$.
 The right-hand side $UV^*$ is a random rank-two matrix, and the target accuracy is $10^{-8}$.
The discretization uses $N^d$ points in the unit square or cube, with equal spacing $h=1/(N+1)$, and we have
\[
  \Lambda(A) \subseteq \bigl[4dh^{-2}\sin^2(\pi h/2),\,4dh^{-2}\cos^2(\pi h/2)\bigr]\approx\bigl[d\pi^2,\,4dh^{-2}\bigr],
  \quad
  \Lambda(B) \subseteq [-10^4,-1].
\]
These are two disjoint intervals, so by \cref{lem:asym_rates} the type-$(p,0)$ bound decays like $e^{-p\mathscr{F}(1)}$, where $\mathscr{F}(1)=\operatorname{acosh}x_F\approx\pi h$ and $x_F=(\lambda_{\max}+\lambda_{\min}+2)/(\lambda_{\max}-\lambda_{\min})$. The target accuracy $10^{-8}$ thus needs $p\approx6(N+1)$ steps.

We compare one-sided g-fADI, with its zeros chosen as the Chebyshev points of the first kind on $E=[\lambda_{\min}(A),\lambda_{\max}(A)]$, against a modified Bartels--Stewart algorithm designed for lopsided equations~\cite{simoncini2016computational}, which diagonalizes $B$ and solves one shifted system with $A$ for each of its $n$ eigenvalues, classical fADI with optimal shifts, the extended Krylov subspace method~\cite{simoncini2007}, and two-pass block Lanczos~\cite{frommer2008}. All solves with $A$ use the sparse direct solver in MATLAB. So that all methods are compared at the same accuracy, g-fADI is run first, and its measured residual becomes the target accuracy of the Krylov methods.

\begin{figure}[t]
\centering
\includegraphics{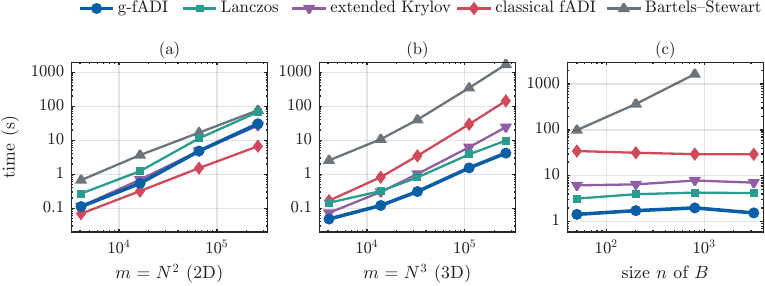}
\caption{The basic test case of Example~1: solve times at the same accuracy, (a) against $m=N^2$ in two dimensions, (b) against $m=N^3$ in three dimensions, and (c) in three dimensions with $N=48$, against the size $n$ of a tridiagonal $B$. ``Lanczos'' is two-pass block Lanczos.}
\label{fig:testbed}
\end{figure}

Three properties decide the outcome. The first is the dimension (\cref{fig:testbed}). On two-dimensional grids, shifted solves with $A$ are cheap, and classical fADI is the fastest method. On three-dimensional grids these solves are more costly, and one-sided g-fADI is the fastest method at every size. %
The other two properties concern $B$. A right step with a tridiagonal $B$ costs $\mathcal{O}(n)$, which is negligible next to a product with $A$. The projection methods, on the other hand, solve, after $k$ steps, a projected Sylvester equation of size $(k\rho)\times n$ (\cref{fig:testbed}). Third, when the spectra of $A$ and $B$ are well separated, g-fADI's iteration count $p$ is small by~\eqref{eq:pcount}. In two dimensions with $N=256$, as the spectrum of $L$ is shifted, with its width kept at $10^4$, so that its left end moves from 0 to $10^6$, the number of steps falls from 1540 to 9, and for a left end of $10^4$ or more, the one-sided iteration is faster than every other method. 

% -------------------------------------------------------------------------------------
\subsection{Example 2: Very nonnormal $B$ and a specialized compression scheme}\label{sec:ex:pde}
Next, we consider a three-dimensional PDE where the problem sizing and properties of the operators give g-fADI a major advantage. Consider the steady state of a population that diffuses through a heterogeneous medium and is also structured by a variable $t\in(0,1)$, such as age or size, which drifts and disperses. Its density satisfies
\begin{equation}\label{eq:age}
  -\nabla\cdot\bigl(\kappa(x)\nabla u\bigr)+v\,\partial_tu-\epsilon\,\partial_t^2u+\mu\,u
  =\phi_1(x)\,\psi_1(t)+\phi_2(x)\,\psi_2(t),
\end{equation}
where $(x,t)\in(0,1)^3\times(0,1)$, with homogeneous Dirichlet conditions, a smooth coefficient $1\le\kappa(x)\le10$, $v=\mu=1$, $\epsilon=0.02$, and Gaussians $\phi_i$ and $\psi_i$, so that the source has rank two.
We discretize with second-order centered finite differences on $N^3$ interior points in $x$, with spacing $h=1/(N+1)$, and on $n=400$ interior points in $t$, with spacing $h_t=1/(n+1)$, and let the $j$th column of $X$ hold $u$ at the $j$th point in $t$. This gives
\begin{equation}\label{eq:agesylv}
  AX-XB=UV^*,\qquad A\in\mathbb{R}^{N^3\times N^3},\quad B=-L^T\in\mathbb{R}^{n\times n},
\end{equation}
where $A$, given in~\eqref{eq:ageA}, is the symmetric positive definite 7-point matrix of $-\nabla\cdot(\kappa\nabla)$ with $\kappa$ evaluated at the midpoints between grid points, and $L$, given in~\eqref{eq:ageL}, is the tridiagonal matrix of $v\,\partial_t-\epsilon\,\partial_t^2+\mu$.

Both operators make this problem difficult, but in different ways. A shifted solve with $A-\beta I$ already costs 103~s and 19~GB for $N=96$, using the MATLAB backslash operator. On the other hand, $B$ is tiny, but the problem is drift-dominated in $t$, and the eigenvector matrix of $B$ has condition number $\sim 7\cdot10^{10}$. Shifted solves with $A$ are too expensive, which rules out Bartels--Stewart and classical fADI. Diagonalizing $B$ would reduce~\eqref{eq:agesylv} to $n$ shifted systems that a multishift Krylov method can solve without any shifted solves, as in Example~4, but in this case would cause a loss of 10 digits of accuracy. One-sided g-fADI avoids all of these issues. It only needs products with $A$ and tridiagonal solves with $B$. We show in \cref{prop:onesided} that we can choose the g-fADI shifts just using a bound on the real part of $W(B)$.

\begin{proposition}\label{prop:onesided}
Let $A$ be normal with $\Lambda(A)\subset E=[a,b]$, and suppose that $W(B)\subset\{z:\operatorname{Re}z\le-f\}$ with $a+f>0$. After $p$ right steps with zeros $\alpha_1,\dots,\alpha_p\in E$ and $X_0=0$, we have, in the 2-norm and in the Frobenius norm,
\begin{equation}\label{eq:onesided}
  \frac{\|X_p-X\|}{\|X\|}\;\le\;\frac{\max_{z\in E}|r_p(z)|}{|r_p(-f)|},\qquad
  r_p(z)=\prod_{j=1}^p(z-\alpha_j).
\end{equation}
In particular, for the Chebyshev points of the first kind on $E$ (the roots of $T_p$, mapped to $E$), the bound is $1/T_p(x_F)$ with $x_F=(a+b+2f)/(b-a)$, the smallest possible (see the proof of \cref{lem:asym_rates}).
\end{proposition}

\begin{figure}[t!]
\centering
\includegraphics[width=0.85\textwidth]{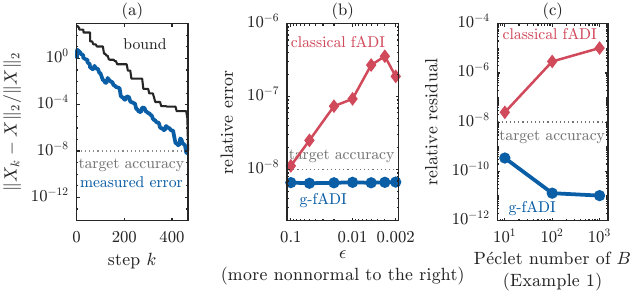}
\caption{(a) The error of one-sided g-fADI for~\eqref{eq:age} with $N=24$ at every step, and the bound~\eqref{eq:onesided}. (b) Increasing the nonnormality of $B$ with $N=32$: the errors of g-fADI and of classical fADI with shifts for $\operatorname{Re}W(B)$, measured against a Schur-based direct solution. The condition number of the eigenvector matrix of $B$ grows from $10^2$ at $\epsilon=0.1$ to $10^{44}$ at $\epsilon=0.005$ and beyond $10^{100}$ at $\epsilon=0.002$. (c) The basic test case of Example~1 in two dimensions ($N=256$) with $L$ replaced by the advection--diffusion matrix~\eqref{eq:tbLadv}: the residuals of the same two methods against the P\'eclet number.}
\label{fig:bound}
\end{figure}

\begin{proof}
By~\eqref{eq:generalized-ADI-error}, $X_p-X=-r_p(A)Xr_p(B)^{-1}$, and $\|r_p(A)\|\le\max_E|r_p|$ since $A$ is normal. For real $\alpha>-f$ and a unit vector $x$, the Cauchy--Schwarz inequality gives $\|(B-\alpha I)x\|\ge|x^*(B-\alpha I)x|\ge \alpha-\operatorname{Re}x^*Bx\ge\alpha+f$. Therefore $\|(B-\alpha_jI)^{-1}\|\le(\alpha_j+f)^{-1}$ for each $j$, and $\|r_p(B)^{-1}\|\le|r_p(-f)|^{-1}$.
\end{proof}
The bound involves no eigenvalues or eigenvectors of $B$ and no constant.
As $B$ is made more nonnormal, until its eigenvector matrix has a condition number beyond $10^{100}$, g-fADI meets the target accuracy every time with essentially the same number of steps (609 to 619), as \cref{prop:onesided} promises. On the other hand, classical fADI with shifts chosen using $\operatorname{Re}W(B)$ misses the target accuracy by up to a factor of 36. The error bound for classical fADI involves all of $W(B)$, which is a more difficult set to approximate (\cref{fig:bound}).
The same happens in the basic test case of Example~1 when $L$ is replaced by the advection--diffusion matrix~\eqref{eq:tbLadv}: g-fADI meets the target accuracy at every P\'eclet number, while classical fADI misses it by up to three orders of magnitude (\cref{fig:bound}).

\begin{figure}[t!]
\centering
\includegraphics{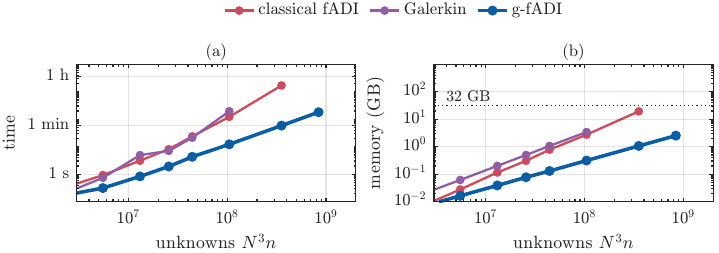}
\caption{(a) Run time and (b) memory held by each method (for its shifted solves, its basis, or for g-fADI the compressed solution and its working arrays) for~\eqref{eq:agesylv} with $n=400$, against the number of unknowns $N^3n$, each method up to the largest size we ran. The dotted line is the 32~GB of the machine.}
\label{fig:cost}
\end{figure}

Since $\operatorname{acosh}x_F\approx2\sqrt{(a+f)/(b-a)}$ and the bound decays like $e^{-p\operatorname{acosh}x_F}$ (\cref{lem:asym_rates}), the number of steps needed for a tolerance $\varepsilon$ is
\begin{equation}\label{eq:pcount}
  p\;\approx\;\tfrac12\sqrt{(b-a)/(a+f)}\;\log(2/\varepsilon).
\end{equation}
For this operator, $a=12h^{-2}\sin^2(\pi h/2)\approx3\pi^2$ and $f \approx 1.2$. The standard Gershgorin bound gives $b\approx120h^{-2}$.  Only $b$ grows with $N$, and we have $\sqrt{(b-a)/(a+f)}\approx1.97\,(N+1)$. For $\varepsilon=10^{-8}$, we get $p\approx18.9\,(N+1)$, independent of $n$ and of how far $B$ is from normal.

\subsubsection{Small-side compression}
With so many steps, memory becomes the bottleneck. 
Each step of g-fADI (Algorithm~\ref{alg:g-fADI}) appends new factors $\hat Z \in \mathbb{C}^{m \times \rho}$ and $\hat Y \in \mathbb{C}^{n \times \rho}$, so after $p$ steps the storage cost is $O(mp\rho + np\rho)$.
If $\rho p > n$, this representation is redundant, since $X \in \mathbb{C}^{m \times n}$ is at most rank-$n$ and is often of much smaller numerical rank.
We compress $X_k$ efficiently by only compressing $Y_k$, which is much smaller than $Z_k$ in this example.

We write $Y_k\approx Q_kG_k$, where $Q_k \in \mathbb{C}^{n\times s_k}$ has orthonormal columns and $G_k \in \mathbb{C}^{s_k \times k\rho}$. $s_k \le \min\{n, k\rho\}$ and often $s_k \ll \min\{n, k\rho\}$.
Then $X_k \approx (Z_kG_k^*)Q_k^*$, and it suffices to store $W_k=Z_kG_k^* \in \mathbb{C}^{m \times s_k}$.
At each step, $\hat Y$ is projected onto the current basis, and its orthogonal component is discarded if it is below a tolerance.
This costs $\mathcal{O}(n s_k\rho)$ per step, and updating $W_k$ costs $\mathcal{O}(m s_k\rho)$. Since $s_k\le n$, neither cost grows with $k$ beyond $\mathcal{O}(n^2\rho)$ and $\mathcal{O}(mn\rho)$.
In practice, $s_k$ is often much smaller than $n$. In this example, $s_k\approx50$ for every $N$, against $n=400$ and $p=2432$ for $N=128$.
Details of small-side compression are given in \cref{alg:onesided}.

\subsubsection{Sampled output}
Often, only a few outputs $\Gamma X$ are wanted, such as the values at a set of points, on a slice, or an average. Since the update is linear, we can track $\Gamma X_k \approx (\Gamma W_k)Q_k^*$ in place of $X_k$, so the memory is $\ell s_k$ numbers for $\ell$ outputs, plus the residual factor $\widetilde Z$ and its product with $A$. Examples~1, 3, and~4 use \cref{alg:onesided} in this form.
We rescale $\widetilde Z$ to unit norm at every step and absorb the scale into $\widetilde Y$, since $r_j(A)U$ otherwise grows beyond the range of double precision within the first hundred steps. Computing this norm is the only operation that combines all $m$ entries of a vector of length $m$.

\begin{algorithm}[h]
\caption{One-sided g-fADI with compression on the small side}\label{alg:onesided}
\KwIn{$A$, $B$, $U$, $V$, zeros $\alpha_1,\dots,\alpha_p$, outputs $\Gamma$ ($\Gamma=I$ for all of $X$), tolerance $\tau$}
\KwOut{$Q$ with orthonormal columns and $W$, where $\Gamma X_p\approx WQ^*$, and a bound $\delta$ on the error of the compression}
$\widetilde Z\leftarrow U$, \ $\widetilde Y\leftarrow V$, \ $Q\leftarrow[\,]$, \ $W\leftarrow[\,]$, \ $\delta\leftarrow0$, \ $\eta\leftarrow0$\;

\For{$j=1,\dots,p$}{
  $\widetilde Y\leftarrow(B^*-\bar\alpha_jI)^{-1}\widetilde Y$ \tcp*{shifted solve with $B^*$}
  $G\leftarrow Q^*\widetilde Y$, \ $Y_\perp\leftarrow\widetilde Y-QG$ \tcp*{part outside the span of $Q$}
  $\eta\leftarrow\max(\eta,\|\Gamma\widetilde Z\|\,\|\widetilde Y\|)$ \tcp*{estimate of $\|\Gamma X_p\|$}
  \eIf{$\|\Gamma\widetilde Z\|\,\|Y_\perp\|>\tau\eta$}{
    $Q_\perp\leftarrow\operatorname{orth}(Y_\perp)$, \ $Q\leftarrow[Q,\ Q_\perp]$, \ $W\leftarrow[W,\ 0]$, \ $G\leftarrow[G^*,\ (Q_\perp^*Y_\perp)^*]^*$\;
  }{
    $\delta\leftarrow\delta+\|\Gamma\widetilde Z\|\,\|Y_\perp\|$
  }
  $W\leftarrow W-(\Gamma\widetilde Z)\,G^*$,\;
  $\widetilde Z\leftarrow(A-\alpha_jI)\,\widetilde Z$ \tcp*{product with $A$}
}
\end{algorithm}

For $N=96$, classical fADI takes 28 minutes and 19~GB, against 1 minute for one-sided g-fADI (\cref{fig:cost}). A Galerkin projection method onto the same spaces must store a basis, already 3.4~GB for $N=64$. Since $A$ is symmetric, two-pass block Lanczos~\cite{frommer2008} avoids storing a basis, but in Example~1 it takes two to four times as long as g-fADI in three dimensions (\cref{fig:testbed}).

% -------------------------------------------------------------------------------------
\subsection{Example 3: 120 coupled PDEs}\label{sec:ex:switch}
In Examples~1 and~2, $A$ is symmetric. Here it is not, and $B$ comes from a Markov chain rather than a discretized derivative. We consider the model in ~\cite{benichou2011}, in which a particle, such as a protein looking for a binding site, moves through $\Omega=(0,1)^3$ under a swirling drift $w$ and switches among $n=120$ internal states. The states are 40 phases of an internal cycle, each with three modes of motion (free, sliding, and nearly immobile) that have their own diffusivity. The probability $u_k(x)$ that a particle starting at $x$ in state $k$ reaches the target $T$, a small ball inside $\Omega$, before it is degraded satisfies
\begin{equation}\label{eq:switch}
  D_k\bigl(\Delta u_k+w\cdot\nabla u_k\bigr)+\sum_{j}R_{kj}\,(u_j-u_k)-\delta\,u_k=0\ \ \text{in }
  \Omega\setminus T,\qquad u_k=1\ \ \text{on }\partial T,
\end{equation}
with no flux through $\partial\Omega$, where $D_k$ is the diffusivity in state $k$, $R$ gives the switching rates, and $\delta>0$ is the degradation rate. We discretize with cell-centered finite differences, second order except at the target, on $N^3$ cells of width $h=1/N$ and put $u_k$ at the cell centers in the $k$th column of $X\in\mathbb{R}^{N^3\times n}$, to obtain
\begin{equation}\label{eq:switchsylv}
  AX\mathcal{D}+XG^T=b\,d^T,\qquad A=H+K,\qquad G=\delta I-R,
\end{equation}
where $\mathcal{D}=\operatorname{diag}(d)$ with $d_k=D_k$, $H$ is the 7-point matrix of $-\Delta$ with the target removed, $K$ is a skew-symmetric discretization of $-w\cdot\nabla$, and $b$ carries the values $u=1$ on $\partial T$ (\cref{app:switch}).
Equation~\eqref{eq:switchsylv} is not immediately amenable to g-fADI, but it can be put in the form of \cref{prop:onesided} by a scaling argument based on the work in~\cite{FortunatoTownsend2020}. The details are given in \cref{app:switch}.

\begin{figure}[t!]
\centering
\includegraphics{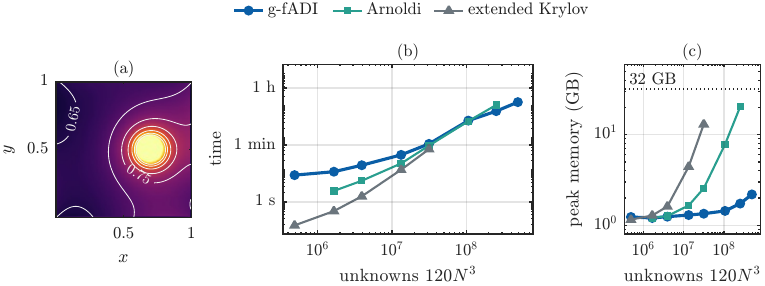}
\caption{(a) The solution $u$ of~\eqref{eq:switch}, averaged over the stationary distribution of the chain, on the layer of cell centers through the target nearest $z=1/2$ for $N=160$. (b) Run time and (c) peak memory against the number of unknowns $120N^3$.}
\label{fig:switch}
\end{figure}

The alternatives are costly here. A shifted solve with $A$ is expensive at these sizes, and Bartels--Stewart would need one for each of the 120 complex eigenvalues of $G$, which has imaginary parts up to $5\cdot10^4$ since the chain is irreversible. Since $A$ is not symmetric, Lanczos is inapplicable, and one has to use Arnoldi, which stores its entire basis. One-sided g-fADI needs only a region that contains $W(A)$ and the bound $f$, and it applies $A$ only through products. We keep only the solution on the slice through the target and its mean over the cube, for every state, so $\Gamma$ in \cref{alg:onesided} has $N^2+1$ rows.
At N=128, the largest grid where Arnoldi also ran, g-fADI is 1.6 times faster and uses 1.7 GB of memory, against 20 GB for Arnoldi.

% -------------------------------------------------------------------------------------
\subsection{Example 4: Collections of shifted linear systems and functions of matrices}\label{sec:ex:shifted}
Many computations require the solutions of a collection of shifted linear systems
\begin{equation}\label{eq:shifted}
  (A+s_jI)x_j=b,\qquad j=1,\dots,k,
\end{equation}
with the same matrix $A\in\mathbb{C}^{m\times m}$ and right-hand side $b$, but different shifts $s_j$. Such systems arise, for instance, when $f(A)b$ is evaluated through a rational approximation $f(z)\approx\gamma_0+\sum_{j=1}^k\gamma_j/(z-\zeta_j)$, computed by the AAA algorithm~\cite{aaa2018} or from a contour integral or quadrature rule~\cite{hale2008, HorningTrefethen2026}. With $s_j=-\zeta_j$, we then have $f(A)b\approx\gamma_0b+\sum_{j=1}^k\gamma_jx_j$. The main alternatives for computing $f(A)b$ are polynomial and rational Krylov methods~\cite{frommer2008, guttel2013, hochbruck1997, saad1992} (see~\cite{higham2008functions} for an overview). When $A$ is large and sparse, solving the systems one at a time requires $k$ shifted solves with $A$. Instead, we introduce  $X\in\mathbb{C}^{m\times k}$ with $x_j$ as its $j$th column, and rewrite~\eqref{eq:shifted} as the Sylvester matrix equation
\begin{equation}\label{eq:shiftsylv}
  AX+XS=b\,\mathbf{1}^T,\qquad S=\operatorname{diag}(s_1,s_2,\dots,s_k),
\end{equation}
which g-fADI can readily be applied to. If $A$ is difficult to invert, we can use a type-$(p,0)$ Zolotarev rational function, i.e., take right steps only, so that the only inverts involved are of the diagonal matrix $S$. The overall cost for $p$ iterations is then $\mathcal{O}(p(M_A+k))$, where $M_A$ is the cost of a matrix-vector product with $A$, and the shifted systems are solved in an inverse-free way. Since $B=-S$ is diagonal, \cref{prop:onesided} applies with $f=\min_j\operatorname{Re}s_j$, and with Chebyshev zeros the number of steps is~\eqref{eq:pcount}, independent of the number of shifts and of their imaginary parts.

Unlike Krylov methods, g-fADI fixes its number of steps in advance from an a priori error bound, so it needs no stopping test. When $A$ is nonnormal, its convergence is governed by the geometry of $W(A)$ through~\eqref{eq:generalized-ADI-error-bound-nonnorm}, not by the eigenvalues of $A$. The method also applies without change when $S$ is not diagonal but, for instance, bidiagonal or banded, as in \cref{sec:intro:teaser} and Example~2.

\subsubsection{Shifted systems from a frequency sweep}\label{sec:ex:sweep}
We now consider a collection of shifted linear systems~\eqref{eq:shifted} with a nonsymmetric $A$, as in Example~3, and shifts $s_j=i\omega_j$ on the imaginary axis. It comes from the frequency response of a heated composite material. For a heat source oscillating at frequency $\omega$, the amplitude $u$ of the steady oscillation satisfies
\begin{equation}\label{eq:sweep}
  -\nabla\cdot\bigl(\kappa(x)\nabla u\bigr)+v\cdot\nabla u+i\omega u=\phi\ \ \text{in }(0,1)^3,\quad u=0\ \ \text{on the boundary},
\end{equation}
where the conductivity $\kappa(x)$ is $10^3$ inside 27 spherical inclusions and $1$ elsewhere, $v$ is a constant flow, and $\phi(x)=\exp(-|x-x_0|^2/0.01)$ is a localized heat source. Suppose we want $u$ at 100 sensors for 1000 frequencies. We discretize~\eqref{eq:sweep} with second-order centered finite differences on $N$ interior points in each direction, using the 7-point stencil for $-\nabla\cdot(\kappa\nabla u)$, with $\kappa$ evaluated midway between adjacent points, and centered differences for $v\cdot\nabla u$. For each frequency $\omega_j$, this gives a shifted linear system $(A+i\omega_jI)x_j=b$, where $x_j$ and $b$ hold the values of $u$ and $\phi$ at the grid points, and $A=H+K$ is given in~\eqref{eq:sweepA} in \cref{app:matrices}. Here, $A$ is sparse but not symmetric, because of the advection term. These 1000 systems are~\eqref{eq:shiftsylv} with $S=\operatorname{diag}(i\omega_1,\dots,i\omega_{1000})$.

\begin{figure}[t!]
\centering
\includegraphics{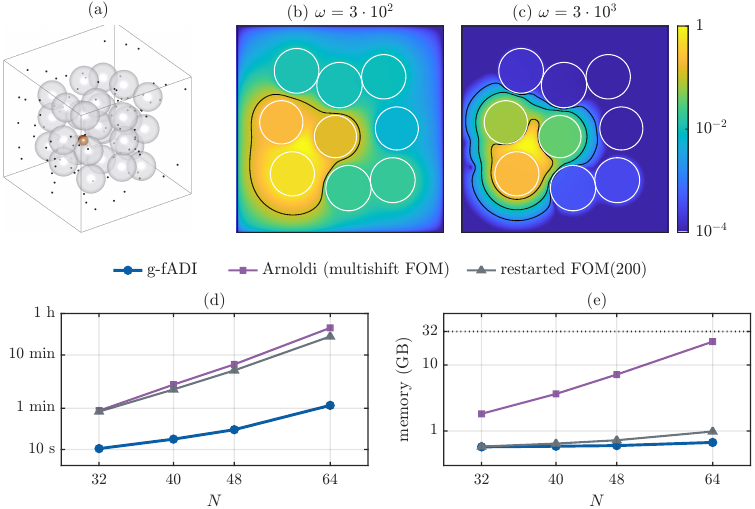}
\caption{(a) The composite material in~\eqref{eq:sweep}: the conducting inclusions (gray), the heat source (orange) and the 100 sensors (black). (b,c) $|u|/\max|u|$ (log scale) on the plane $z=1/2$ through the source at $\omega=3\cdot10^2$ and $\omega=3\cdot10^3$ ($N=127$), with isolines at $10^{-1}$ and, in (c), $10^{-2}$ and $10^{-3}$. The white circles are the inclusions cut by the plane. (d) Run time and (e) peak memory for the response at the sensors at 1000 frequencies, against the grid size $N$, for one-sided g-fADI, for Arnoldi with one basis for all frequencies (the multishift full orthogonalization method), and for the same method restarted every 200 steps, all to a relative residual of $10^{-8}$. All setup costs for g-fADI are included in the timings.}
\label{fig:sweep}
\end{figure}

For $N=64$ points in each direction, one-sided g-fADI, computing only the sensor values rather than the whole solution (\cref{alg:onesided} with $\Gamma$ the 100 sensor rows), computes all $10^5$ of them in 68~s and 0.7~GB. Arnoldi with one basis for all frequencies takes 32 minutes and 22.5~GB, and restarting it keeps the memory small but still takes 22 minutes (\cref{fig:sweep}).

In this example, g-fADI needs four to six times as many products with $A$ as Arnoldi, since its shifts are fixed in advance rather than adapted to $A$ and $b$, but every step costs the same. In contrast, every new Arnoldi vector must be orthogonalized against, and stored with, all of the previous ones. Arnoldi needs more steps as $N$ grows, so its cost grows quadratically in the number of steps.

% -------------------------------------------------------------------------------------

\section*{Acknowledgments and AI disclosure statement} We are appreciative of Nick Trefethen, Yuji Nakatsukasa, and several others in attendance at the BIRS conference ``Challenges, Opportunities, and New Horizons in Rational Approximation''  for several conversations related to generalized Zolotarev rationals.  We used Claude Opus 5.5 in various ways: the PDE-based examples and applications were initially suggested by Claude and then refined by us. \Cref{prop:onesided,prop:scaling} were suggested, along with their proofs, entirely by Claude, and then rigorously checked and tested by the authors. Claude Code was also fully integrated into our implementations of the algorithms, numerical testing, and the creation of our GitHub repository~\cite{gfadiCode}, where all our code and examples can be found. As authors we assume full responsibility for all results and code contained in this manuscript. 
% =====================================================================================

\appendix
\section{The matrices in the experiments}\label{app:matrices}
All grids consist of equispaced points, and discretizations in $d$ dimensions are unfolded lexicographically via Kronecker products. We write $\operatorname{tridiag}(\alpha,\beta,\gamma)$ for the Toeplitz matrix with $\alpha$, $\beta$, and $\gamma$ as the subdiagonal, diagonal, and superdiagonal, respectively.

\subsection{The Laplacian}
With $N$ interior points per direction having spacing $h=1/(N+1)$, the second-order centered finite-difference matrix for $-\Delta$ on $(0,1)^d$ with homogeneous Dirichlet conditions is
\begin{equation}\label{eq:lap}
  A=I\otimes T+T\otimes I\quad(d=2),\qquad
  A=I\otimes I\otimes T+I\otimes T\otimes I+T\otimes I\otimes I\quad(d=3),
\end{equation}
where $T=h^{-2}\operatorname{tridiag}(-1,2,-1)\in\mathbb{R}^{N\times N}$. It is symmetric positive definite with eigenvalues in $[4dh^{-2}\sin^2(\pi h/2),\,4dh^{-2}\cos^2(\pi h/2)]$.

\subsection{The basic test case and the model problem}
In Example~1, $A$ is~\eqref{eq:lap} and $B=-L^T$ with
\begin{equation}\label{eq:tbL}
  L=I+\tfrac{10^4-1}{4}\operatorname{tridiag}(-1,2,-1)\in\mathbb{R}^{200\times200},
\end{equation}
the second-order finite-difference matrix for $1-\epsilon\,\partial_t^2$ on $(0,1)$ with homogeneous Dirichlet conditions and $\epsilon=(10^4-1)/(4\cdot201^2)$. Its spectrum lies in $[1,10^4]$. In \cref{fig:bound}, $L$ is replaced by the nonnormal matrix
\begin{equation}\label{eq:tbLadv}
  L=I+\sigma\bigl(h_n^{-2}\operatorname{tridiag}(-1,2,-1)+v\,h_n^{-1}\operatorname{tridiag}(-1,1,0)\bigr),\qquad h_n=1/201,
\end{equation}
with second-order centered differences for the diffusion term, first-order backward differences for the advection term, P\'eclet number $v/2$, and $\sigma$ chosen so that the Gershgorin bound on the spectrum is $10^4$.

In the model problem of \cref{sec:intro:teaser}, $A=H+K$, where $H$ is~\eqref{eq:lap} and
\begin{equation}\label{eq:tbK}
  K=w_1\,I\otimes I\otimes D_{\mathrm c}+w_2\,I\otimes D_{\mathrm c}\otimes I+w_3\,D_{\mathrm c}\otimes I\otimes I,
\end{equation}
with $D_{\mathrm c}=(2h)^{-1}\operatorname{tridiag}(-1,0,1)$, is the second-order centered finite-difference matrix for $w\cdot\nabla$ with $w=5\,(1,0.6,0.3)/\|(1,0.6,0.3)\|$. It is skew-symmetric.

\subsection{The structured population (Example~2)}
We use $N$ interior points per direction in $x$ with $h=1/(N+1)$, and $n=400$ interior points in $t$ with $h_t=1/(n+1)$. The $j$th column of $X$ holds $u$ at the $j$th point in $t$. In~\eqref{eq:agesylv},
\begin{equation}\label{eq:ageA}
  A=\sum_{d=1}^3D_d^T\operatorname{diag}(\kappa_d)\,D_d,
\end{equation}
where $D_1=I\otimes I\otimes D$, $D_2=I\otimes D\otimes I$, and $D_3=D\otimes I\otimes I$. Here $D\in\mathbb{R}^{(N+1)\times N}$ has $h^{-1}$ on its diagonal and $-h^{-1}$ on its subdiagonal, so that $(Du)_j=(u_j-u_{j-1})/h$, with $u_0=u_{N+1}=0$, approximates $u'$ at the midpoint between the $(j-1)$st and $j$th grid points, and $\kappa_d$ holds the values of $\kappa$ at the midpoints in direction $d$. This is the symmetric positive definite second-order 7-point matrix of $-\nabla\cdot(\kappa\nabla)$. The small coefficient is $B=-L^T$ with
\begin{equation}\label{eq:ageL}
  L=\epsilon\,h_t^{-2}\operatorname{tridiag}(-1,2,-1)+v\,(2h_t)^{-1}\operatorname{tridiag}(-1,0,1)+\mu I\in\mathbb{R}^{n\times n},
\end{equation}
the second-order centered finite-difference matrix for $v\,\partial_t-\epsilon\,\partial_t^2+\mu$.

\subsection{The switching diffusion (Example~3)}\label{app:switch}
We use a discretization with~$N^3$ cells of width $h=1/N$. With no flux through the walls, the second-order cell-centered finite-difference matrix for $-\Delta$ is
\begin{equation}\label{eq:switchlap}
  I\otimes I\otimes T_{\mathrm N}+I\otimes T_{\mathrm N}\otimes I+T_{\mathrm N}\otimes I\otimes I,
\end{equation}
where $T_{\mathrm N}$ is $h^{-2}\operatorname{tridiag}(-1,2,-1)$ with its first and last diagonal entries replaced by $h^{-2}$. To obtain $H$, we remove the cells whose centers lie in the target and impose $u=1$ at those centers, which moves these values into $b$. This is first order at the target, since it ignores the offset between the cell centers and $\partial T$. The drift term is discretized by second-order centered differences, with $w$ evaluated at the cell faces,
\begin{equation}\label{eq:switchK}
  (Ku)_{\mathbf i}=-(2h)^{-1}\sum_{d=1}^3\Bigl(w_d(x_{\mathbf i+e_d/2})\,u_{\mathbf i+e_d}-w_d(x_{\mathbf i-e_d/2})\,u_{\mathbf i-e_d}\Bigr),
\end{equation}
where $\mathbf i$ indexes the cells, $e_d$ is the $d$th unit vector, and $x_{\mathbf i\pm e_d/2}$ is the midpoint of the face between cell $\mathbf i$ and its neighbor in direction $\pm e_d$. Couplings across the walls and into the target are omitted. Since both cells of a face use the same value of $w_d$, $K$ is skew-symmetric. We take $w=\bigl(-(y-\tfrac12),\,x-\tfrac12,\,0.3\sin2\pi x\bigr)/\sqrt{0.59}$, so that $|w|\le1$ and $\nabla\cdot w=0$, and~\eqref{eq:switchK} is then a second-order discretization of $-w\cdot\nabla u$.

\paragraph{The field of values of $A$}
With $c^2=\sum_d\max|w_d|^2$, the Cauchy--Schwarz inequality gives $|x^*Kx|\le c\,\|x\|\,(x^*Hx)^{1/2}$, and $|x^*Kx|\le\|K\|\,\|x\|^2\le y_{\max}\|x\|^2$ with $y_{\max}=h^{-1}\sum_d\max|w_d|$. Hence
\[
  W(A)\subseteq\Omega_A=\bigl\{z:\ \operatorname{Re}z\in[\lambda_{\min}(H),\lambda_{\max}(H)],\ \
  |\operatorname{Im}z|\le\min\bigl(c\sqrt{\operatorname{Re}z},\,y_{\max}\bigr)\bigr\}.
\]
By Crouzeix's bound~\cite{jin2026numerical,lorist2026}, $\|r_p(A)\|\le2\max_{\Omega_A}|r_p|$, so \cref{prop:onesided} holds with $\Omega_A$ in place of $E$ and an extra factor of~$2$. We take the zeros to be the Chebyshev points of the first kind on an interval $[a',b']$ with $b'\ge\lambda_{\max}(H)$, and choose $a'$, $b'$, and $p$ by a scalar search so that $2\max_{\Omega_A}|r_p|/|r_p(-f)|\le\varepsilon=10^{-8}$ with $p$ as small as possible.

\paragraph{The scaling of $G$}
The matrix $G=\delta I-R$ is not of the form required by \cref{prop:onesided}. Since the generator $R$ has nonnegative off-diagonal entries and zero row sums, and $\delta>0$, $G$ is a nonsingular $M$-matrix.\footnote{A real matrix with nonpositive off-diagonal entries is a nonsingular $M$-matrix if it is invertible and its inverse is entrywise nonnegative. Every strictly diagonally dominant matrix with positive diagonal and nonpositive off-diagonal entries is one~\cite{berman1994}.} For such matrices, a diagonal similarity moves the field of values into the right half-plane~\cite{berman1994}. We use the following explicit version.

\begin{proposition}\label{prop:scaling}
Let $G\in\mathbb{R}^{n\times n}$ be a nonsingular $M$-matrix, let $x=G^{-1}\mathbf{1}$ and $y=G^{-T}\mathbf{1}$, and let $S=\operatorname{diag}(\sqrt{y_i/x_i})$. Then $L=SGS^{-1}$ satisfies
\[
  \operatorname{Re}W(L)\;\ge\;f:=\tfrac12\min_i\bigl(1/x_i+1/y_i\bigr)>0 .
\]
\end{proposition}

\begin{proof}
Since $G^{-1}$ is nonnegative and nonsingular, $x,y>0$. The symmetric matrix $J=L+L^T$ has nonpositive off-diagonal entries, and $z=Sx>0$ satisfies $Jz=S^{-1}(S^2Gx+G^TS^2x)=S^{-1}(S^2\mathbf{1}+G^Ty)$, so that $(Jz)_i/z_i=1/x_i+1/y_i$. The Collatz--Wielandt formula~\cite{berman1994}, applied to the nonnegative matrix $\nu I-J$ for large $\nu$, gives $\lambda_{\min}(J)\ge\min_i(Jz)_i/z_i$, and the claim follows from $\operatorname{Re}W(L)\ge\tfrac12\lambda_{\min}(J)$.
\end{proof}

The matrix $G'=\mathcal{D}^{-1/2}G\mathcal{D}^{-1/2}$ is again a nonsingular $M$-matrix. With its scaling $S$, the substitution $X=\widetilde X\mathcal{D}^{-1/2}S^{-1}$ turns~\eqref{eq:switchsylv} into $A\widetilde X+\widetilde XL^T=bv^T$ with $L=SG'S^{-1}$ and $v=S\mathcal{D}^{-1/2}d$, and \cref{prop:scaling} gives $f=0.22$.

\subsection{The frequency sweep (Example~4)}
We use $N$ interior points per direction with $h=1/(N+1)$, and $A=H+K$ with
\begin{equation}\label{eq:sweepA}
\begin{aligned}
  H&=\sum_{d=1}^3D_d^T\operatorname{diag}(\kappa_d)\,D_d,\\
  K&=v_1\,I\otimes I\otimes D_{\mathrm c}+v_2\,I\otimes D_{\mathrm c}\otimes I+v_3\,D_{\mathrm c}\otimes I\otimes I,
\end{aligned}
\end{equation}
where $D_d$ and $\kappa_d$ are as in~\eqref{eq:ageA}, $D_{\mathrm c}$ is as in~\eqref{eq:tbK}, and $v=(1,0.6,0.3)/\|(1,0.6,0.3)\|$. Here $H$ is the symmetric positive definite second-order 7-point matrix of $-\nabla\cdot(\kappa\nabla)$, and $K$ is the skew-symmetric second-order centered finite-difference matrix for $v\cdot\nabla$.

\section{The type-\texorpdfstring{$(p,q)$}{(p,q)} Z3 algorithm}\label{app:z3alg}
\Cref{alg:modified aaa} gives the AAA part of the method in \cref{sec:computing shifts}. It approximates $f=\sgn_{E/F}$ by a rational function $g$ in barycentric form whose transform $\tilde S[g]$ has type $(p,q)$.
\begin{algorithm}[H]
\caption{Type-$(p,q)$ Z3 algorithm}
\label{alg:modified aaa}
\begin{algorithmic}[1]
\Require sample points $Z=\{z_1,\dots,z_M\}\subset E\cup F$, $f_i=f(z_i)$ with $f=-1$ on $E$ and $f=1$ on $F$, type $(p,q)$, tolerance $\mathrm{tol}$
\Ensure Rational function $r \in \RR_{p,q}$

\State $S\gets\emptyset$, \ $g \gets \text{average}(f)$
\For{$k = 1,2,\dots,\max{\{p,q\}}+1$}
    \State $j \gets \arg\max_{i\not\in S} |f_i - g(z_i)|$, \ $S \gets S \cup \{j\}$
    \State Support points $\{s_i\}_{i=1}^k = \{z_j : j\in S\}$, other points $\{\zeta_i\}_{i=1}^{M-k} = \{z_i : i \not\in S\}$
    \State Set
        \begin{equation*}
          g(z) \;=\; \Big(\sum_{i=1}^k \frac{w_i f(s_i)}{z-s_i}\Big)
                \bigg/ \Big(\sum_{i=1}^k \frac{w_i}{z-s_i}\Big),
          \qquad
          w=\Pi\hat w,\quad \hat w = \underset{\|\Pi\hat w\| = 1}{\textup{argmin}}\; \|\mathbb{L}\Pi\hat w\|,
        \end{equation*}
        where $\mathbb{L}_{i,j} = \big(f(s_j)-f(\zeta_i)\big)/(\zeta_i-s_j)$ and $\Pi$ is the projector onto the nullspace of $P$ in~\eqref{eq:P-matrix}
    \State \textbf{if} $\max_{i} |f_i - g(z_i)|\le \mathrm{tol}$ \textbf{then break}
\EndFor
\State \Return $r = \tilde S[g] = (1+g)/(1-g)$
\end{algorithmic}
\end{algorithm}
\noindent The constraint matrix $P\in \C^{|p-q| \times k}$ makes the leading coefficients of the numerator or the denominator of $\tilde S[g]$ vanish. Its entries are
    \begin{equation}\label{eq:P-matrix}
        P_{i,j} = \big(f(s_j) + \sgn(q-p)\big)\cdot c_{k-i,j},
    \end{equation}
    where the $c_{\nu,j}$ are the coefficients of
    \begin{equation}\label{eq:ellz}
        \ell_j(z)
        = \prod_{\ell\neq j}(z-s_\ell)
        = \sum_{\nu=0}^{k-1}c_{\nu,j} z^\nu,
        \qquad 1 \le j \le k.
    \end{equation}

\bibliographystyle{siamplain}
\bibliography{refs}

\end{document}